\documentstyle[12pt]{article}
\newcommand{\la}[1]{\mbox{\large $#1$}}
\newcommand{\La}[1]{\mbox{\Large $#1$}}
\newcommand{\LA}[1]{\mbox{\LARGE $#1$}}

\newcommand{\mb}[1]{\mbox{#1}}

\newcommand{\lnab}[1]{\la{\nabla}_{\!\!#1}} 
\newcommand{\hlnab}[1]{\hat{\la{\nabla}}_{\!\!#1}} 
\newcommand{\lnabp}[1]{\la{\nabla'}_{\!\!#1}} 
\newcommand{\lnabo}[1]{\la{\nabla}^{\bot}_{\!\!#1}} 

\newcommand{\qed}{\mbox{~~~\boldmath $\Box$}} 
\newcommand{\R}[1]{I\!\!R^{#1}}
\newcommand{\Co}{l\!\!\! C}
\newcommand{\ra}{\rightarrow}
\newcommand{\non}{\nonumber}

\newcommand{\al}{\alpha}
\newcommand{\be}{\beta}
\newcommand{\ga}{\gamma}
\newcommand{\bal}{\bar{\alpha}}
\newcommand{\bbe}{\bar{\beta}}
\newcommand{\bga}{\bar{\gamma}}
\newcommand{\bmu}{\bar{\mu}}
\newcommand{\brho}{\bar{\rho}}
\newcommand{\tg}{\tilde{g}}
\newcommand{\gdf}[3]{g(\mbox{\large ${\nabla}$}_{\!\!{#1}} dF({#2}),JdF({#3}))}
\newcommand{\gf}[3]{g_{#1}{#2}{#3}}
\newcommand{\Jw}{J_{\omega}}
\newcommand{\kw}{{\cal K}_{\omega}}

\newcommand{\Fw}{F^{*}\omega}
\newtheorem{Lm}{Lemma}
\newtheorem{Pp}{Proposition}
\newtheorem{Th}{Theorem}
\newtheorem{Cr}{Corollary}

\begin{document}
\baselineskip .6cm
\setcounter{page}{1}
\title{Minimal submanifolds of K\"{a}hler-Einstein  manifolds with equal K\"{a}hler angles } \author{ Isabel M.\ C.\ Salavessa\raisebox{.8ex}{\scriptsize 1} 
 and Giorgio Valli\raisebox{.8ex}{\scriptsize \dag2} }
\protect\footnotetext{Deceased on October 2nd, 1999}
\date{}
\maketitle
\mbox{ } \\[-15mm]
 {\footnotesize 1 Centro de F\'{\i}sica das
Interac\c{c}\~{o}es Fundamentais, Instituto Superior T\'{e}cnico,
Edif\'{\i}cio Ci\^{e}ncia,\\[-2mm] Piso 3,
1049-001 LISBOA, Portugal;~~
e-mail: isabel@cartan.ist.utl.pt}\\[-1mm]
{\footnotesize 2 Dipartimento di
Matematica, Universit\`{a} di Pavia,
Via Abbiategrasso 215,
27100 PAVIA, Italy}\\[1cm]
{\small {\bf Abstract:} We consider $F: M\ra N $ a minimal oriented compact real 
$2n$-submanifold $M$, immersed
into a K\"{a}hler-Einstein manifold $N$ of 
complex dimension $2n$, and scalar curvature $R$.
We assume that $n\geq 2$ and  $F$ has equal K\"{a}hler 
angles. Our main result is to prove that, if $n=2$ and $R\neq 0$, then
 $F$ is either a complex submanifold or a Lagrangian submanifold. 
We also prove
that, if $n\geq 3$ and $F$ has  no complex points, then :
(A) If  $R<0$,  then $F$ is Lagrangian; 
(B) If $R=0$, the  K\"{a}hler angle must be constant. 
We also study pluriminimal submanifolds with equal K\"{a}hler angles,
and prove that, if they are not complex submanifolds, $N$ must be  Ricci-flat and there is a natural parallel homothetic isomorphism between $TM$ and  the normal bundle.\\[1mm] 
}
\markright{\sl\hfill  Salavessa - Valli \hfill} 
\section{Introduction}
\setcounter{Th}{0}
\setcounter{Pp}{0}
\setcounter{Cr} {0}
\setcounter{Lm} {0}
\setcounter{Def} {0}
\setcounter{equation} {0}
Let $(N,J,g)$ be a K\"{a}hler manifold of complex dimension $2n$ and 
 $F:M\ra N$ an immersed submanifold of real dimension $2n$. 
We denote by $\omega$ the K\"{a}hler form of $N$, $\omega(X,Y)=
g(JX,Y)$. On $M$ we take the induced metric $g_{M}= F^{*}g$.  $N$ is 
K\"{a}hler-Einstein  if its  Ricci tensor is a multiple of
the metric, $Ricci^{N}= Rg$. 
 At each point $p\in M$,  we identify
$F^{*}\omega$ with a skew-symmetric operator of $T_{p}M$ by using the 
musical isomorphism with respect to $g_{M}$, namely
$g_{M}(F^{*}\omega(X),Y)=F^{*}\omega(X,Y)$. We take its polar 
decomposition
\begin{equation}
F^{*}\omega =  \tg \Jw
\end{equation}
where $J_{\omega}:T_{p}M\ra T_{p}M$ is a ( in fact unique) partial isometry with the same 
kernel ${\cal K}_{\omega}$ as  of $F^{*}w$, and where $\tg$ is the 
positive semidefinite 
operator $\tg = |F^{*}\omega|=\sqrt{-(F^{*}\omega)^{2} }$. 
It turns out that
$J_{\omega}:{\cal K}^{\bot}_{\omega}\ra{\cal K}^{\bot}_{\omega}$ 
defines a complex structure  on ${\cal K}^{\bot}_{\omega}$, the
orthogonal compliment of $\kw$ in $T_{p}M$. Moreover, it is
$g_{M}$-orthogonal. 
If we denote by $\Omega^{0}_{2k}$ the largest open set of
 $M$ where $F^{*}\omega$ has constant rank $2k$, 
$0\leq k\leq n$, then ${\cal K}^{\bot}_{\omega}$ is a smooth sub-vector 
bundle of $TM$ on $\Omega^{0}_{2k}$. Moreover,  $\tg$ and $J_{\omega}$ 
are both smooth on these open sets. The tensor
$\tg$ is continuous on all $M$ and locally Lipschitz, for
the map $P\ra |P|$ is Lipschitz in the space of normal operators. 
Let $\{ X_{\al},Y_{\al}\}_{1\leq \al \leq n}$ be a
 $g_{M}$-orthonormal  basis of $T_{p}M$, that diagonalizes 
$F^{*}\omega$ at $p$, that is 
\begin{equation}
F^{*}\omega=\bigoplus_{1\leq \al\leq n}  \left[ \begin{array}{cc}
             0 &  -\cos\theta_{\al} \\
             \cos\theta_{\al} &  0 \end{array} \right], 
\end{equation}
where $\cos\theta_{1}\geq \cos\theta_{2}\geq \ldots \geq \cos\theta_{n}
\geq 0$. The angles  $\{\theta_{\al}\}_{ 1\leq \al\leq n}$
 are the \em K\"{a}hler angles  \em of $F$ at $p$. Thus, $\forall \al$, $F^{*}\omega(X_{\al})=\cos\theta_{\al}Y_{\al}$,
$F^{*}\omega(Y_{\al})=-\cos\theta_{\al}X_{\al}$ 
and if $k\geq 1$,
where $2k$ is the rank of $F^{*}\omega$ at $p$,
$J_{\omega}X_{\al}=Y_{\al}$ $~\forall \al \leq k$. The Weyl's
perturbation theorem applied to the eigenvalues of the symmetric 
operator $|\Fw |$ shows that, ordering the $\cos\theta_{\al}$ in the
above way, the map 
$p\ra cos\theta_{\al}(p)$ is locally Lipschitz on $M$, for each $\al$.
A \em complex direction \em  of $F$ is a real two-plane $P$ of $T_{p}M$
such that $dF(P)$ is a complex line of $T_{F(p)}N$, i.e., 
$JdF(P)\subset dF(P)$. Similarly, $P$ is said to be a \em Lagrangian 
direction \em  of $F$ if $\omega$ va\-ni\-shes on $dF(P)$,
that is, $JdF(P)\bot dF(P)$. The immersion $F$ has no complex directions  iff
$\cos\theta_{\al}<1$ $\forall \al$. 
 $M$ is a complex submanifold
iff  $cos\theta_{\al}=1$ $\forall \al$,  and is a 
Lagrangian submanifold iff  $cos\theta_{\al}=0$ $\forall \al$.
We say that $F$ has \em equal K\"{a}hler angles \em if $\theta_{\al}
=\theta$ $\forall \al$. Complex and Lagrangian submanifolds are
examples of such case. If $F$ is a complex submanifold, then $\Jw$ is the
complex structure induced by $J$ of $N$.
The K\"{a}hler angles  are some functions
that at each point $p$ of $M$ measure the deviation of the tangent plane
$T_{p}M$ of $M$ from a complex or a Lagrangian subspace of $T_{F(p)}N$.
This concept was introduced by Chern and Wolfson [Ch-W] for surfaces,
namely $\Fw=\cos\theta\, Vol_{M}$. This $\cos\theta$ may have negative values
and is smooth on all $M$. In our definition, for $n=1$, we demanded
$\cos\theta\geq 0$, that is, it is the modulus of the $\cos\theta$
given for surfaces. This may make our  $\cos\theta$ do
 not  be smooth. We have chosen this definition, because in higher
dimensions we do not have a preferential orientation assigned to the real
planes $span\{X_{\al}, Y_{\al}\}$.\\[-2mm]

Our main aim is  to find conditions for a minimal submanifold $F$
to be Lagrangian or complex, or $M$ to be a K\"{a}hler manifold with respect
to $\Jw$. The first result in this direction is due to 
Wolfson, for the case $n=1$:
\begin{Th} {\bf [W]} 
 If  $M$ is a real compact surface and $N$ is a complex
K\"{a}hler-Einstein surface with $R<0$, anf if $F$ is minimal
with no complex points, then $F$ is Lagrangian.
\end{Th}
Some results of [S-V]  are a generalization of the above theorem to 
higher dimensions. In this paper we study the case of 
 equal K\"{a}hler angles.  Let us denote by $\lnab{X}dF(Y)=\lnab{}dF(X,Y)$
the second fundamental form of $F$. It  is a symmetric tensor and 
takes values in the normal bundle
 $NM=(dF(TM))^{\bot}$. $F$ is minimal iff
$trace_{g_{M}}\lnab{}dF =0$. Let $(~~)^{\bot}$ denote
 the orthogonal projection of $F^{-1}TN$ onto the normal bundle.
If $F$ is an immersion with no complex directions  at $p$ 
and $\{X_{\al},Y_{\al}\}$ diagonalizes $F^{*}\omega$
at $p$, then $\{dF(Z_{\al})$, $
dF(Z_{\bal})$, $(JdF(Z_{\al}))^{\bot}$, $(JdF(Z_{\bal}))^{\bot}\}$ 
constitutes a complex basis of $T^{c}_{F(p)}N$, where \\[-3mm]
\begin{equation}
Z_{\al}=\frac{X_{\al}-iY_{\al}}{2}=``\al ", ~~~~~~~~~Z_{\bal}=\overline{Z_{\al}}=\frac{X_{\al}+iY_{\al}}{2}=
``\bal "
\end{equation}
are complex vectors of the complexified tangent space of $M$
at $p$. We extend to the complexified vector bundles the 
Riemannian tensor metric
 $g_{M}$ (sometimes denoted by $\langle,\rangle$), the curvature tensors
of $M$ and $N$, and any other tensors that will occur, always
by $\Co$-multilinearity.
 On $M$ the Ricci tensor of $N$ can be described by the following 
expression ([S-V]): for $U,V\in T_{F(p)}N$,\\[-3mm]
\begin{equation}
Ricci^{N}(U,V)=\sum_{1\leq \mu\leq n}\frac{4}{\sin^{2} \theta_{\mu}}
R^{N}(U,JV,dF(\mu),(JdF(\bmu))^{\bot}),
\end{equation}
where $R^{N}$ denotes the Riemannian curvature tensor of $N$.
An application of Codazzi equation to the above 
expression proves that, if $N$ is K\"{a}hler-Einstein with $R\neq 0$,  
Theorem 1.1 can be generalized to any  dimension 
for totally geodesic maps ([S-V]). \\[-2mm]

We can also obtain the same conclusion to ``broadly-pluriminimal" 
immersions  for $n=2$,
and $N$ K\"{a}hler-Einstein  with negative Ricci tensor ([S-V]). 
A minimal immersion $F$ is said to be \em broadly-pluriminimal, \em 
if, for each $p\in\Omega_{2k}^{0}$, with $ k\geq 1 $, 
$F$ is pluriharmonic with respect to any $g_{M}$-orthogonal
complex structure  $\tilde{J}=\Jw\oplus J'$ on $T_{p}M$ 
 where $J'$ is any  $g_{M}$-orthogonal complex structure
of $\kw$ at $p$, that is, $(\lnab{}dF)^{(1,1)}=0$.  
The (1,1)-part of $\lnab{}dF$ is just given by
$(\lnab{}dF)^{(1,1)}(X,Y)=\frac{1}{2}
\La{(}\lnab{}dF(X,Y)+\lnab{}dF(\tilde{J}X,
\tilde{J}Y) \La{)}~~\forall X,Y\in T_{p}M.$
If ${\cal K}_{\omega}=0$, this means that $F$ is pluriharmonic 
with respect to the almost complex structure $J_{\omega}$
(see for example [O-V]). In this case, we say that
$F$ is pluriminimal in the usual sense, or simply \em pluriminimal. \em
Pluriharmonic immersions are obviously minimal.
If $F$ has equal K\"{a}hler angles, then only $\Omega^{0}_{2n}$ is 
considered, where ${\cal K}_{\omega}=0$ and $ \tilde{J}=\Jw$.
Products of minimal
real surfaces of K\"{a}hler surfaces, totally geodesic submanifolds,
 minimal Lagrangian submanifolds, 
and complex submanifolds are examples of broadly-pluriminimal submanifolds.
We will see in sections 2 and 3 that the concept of 
broadly-pluriminimality, for immersions without complex directions
and with equal K\"{a}hler angles, may 
have a geo\-me\-tric interpretation in  terms of the torsion of a new 
Riemannian connection on $TM$,
described through an isomorphism $\Phi$ from the tangent bundle of
$M$ into the normal bundle.  Pluriminimal immersions with equal K\"{a}hler angles immersed into K\"{a}hler-Einstein manifolds,
that are not complex submanifolds, have constant K\"{a}hler angle, and only exist on Ricci-flat manifolds.  In this case, $\Phi$ defines a parallel
homothetic isomorphism between $TM$ and $NM$.\\[-2mm] 

For a minimal immersion $F$ with no complex directions  we consider
the locally Lipschitz map, symmetric on the K\"{a}hler angles,\\[-3mm]
\begin{equation}
\kappa=\sum_{1\leq \al\leq n}log\left(\frac{1+\cos\theta_{\al}}
{1-\cos\theta_{\al}}\right).
\end{equation}
This map is smooth on each $\Omega_{2k}^{0}$, non-negative,
and vanishes at Lagrangian points. It is an increasing map on each 
$\cos\theta_{\al}$. In [S-V] we have given an
 expression for $\triangle \kappa$  at a point $p_{0}\in \Omega_{2k}^{0}$,  which we prove in the appendix of this paper, namely,\\[-7mm]
\begin{eqnarray}
\triangle \kappa
&=&4i\sum_{\be} Ricci^{N}(JdF({\be}),dF({\bbe}))\\[-2mm]
&&+\sum_{\be,\mu}\!\frac{32}{\sin^{2}\theta_{\mu}}Im \La{(}\!
R^{N}(dF({\be}),dF({\mu}),dF({\bbe}),
JdF({\bmu})\!+\!i\cos\theta_{\mu}dF(\bmu))\!\La{)}\non\\[-2mm]
&&-\sum_{\be,\mu,\rho}\!\!\!\frac{64(\cos\theta_{\mu}
\!+\!\cos\theta_{\rho})}{\sin^{2}\theta_{\mu}\sin^{2}\theta_{\rho}}
Re\la{(}\gdf{\be}{\mu}{\brho}\gdf{\bbe}{\rho}{\bmu}\!\la{)}
\nonumber\\[-2mm]
&& +\sum_{\be,\mu,\rho}
\frac{32(\cos\theta_{\rho}-\cos\theta_{\mu})}
{\sin^{2}\theta_{\mu}\sin^{2}\theta_{\rho}}\;(|\gdf{\be}{\mu}{\rho}|^{2}
+|\gdf{\bbe}{\mu}{\rho}|^{2})\non\\[-1mm]
&&+\sum_{\be,\mu,\rho}\frac{32(\cos\theta_{\mu}+\cos\theta_{\rho})}
{\sin^{2}\theta_{\mu}}\,\LA{(}
|\langle\lnab{{\be}}{\mu},{\rho}\rangle|^{2} +
|\langle\lnab{{\bbe}}{\mu},{\rho}\rangle|^{2}\LA{)},\nonumber 
\end{eqnarray}
where $\{X_{\al}, Y_{\al}\}_{1\leq \al\leq n}$ is a $g_{M}$-orthonormal
local frame of $M$,  with $Y_{\al}=J_{\omega}X_{\al}$ for $\al\leq k$,
$\{X_{\al},Y_{\al}\}_{\al\geq k+1}$ any $g_{M}$-orthonormal frame of ${\cal K}_{\omega}$, and which at $p_{0}$ diagonalizes $F^{*}\omega$.
For $F$ pluriminimal on  $\Omega_{2n}^{0}$ and $N$ K\"{a}hler-Einstein , 
we can get the following very simple final expression on $\Omega^{0}_{2n}$
([S-V]) 
\begin{equation}
\triangle\kappa=-2R \La{(}\sum_{1\leq\be\leq n}\cos\theta_{\be}\La{)}.
\end{equation}
If $F$ has equal K\"{a}hler angles, then the expression of 
$\triangle\kappa$ given in
(1.6) can also be substantially simplified.  Minimal 
surfaces with constant curvature and constant K\"{a}hler angle in complex
space forms have been classified in [O]. Conditions on the 
curvature of $M$, $N$,  and/or constant  equal K\"{a}hler angles  lead
to some conclusions in our case as well, as we show in the theorems 
below. Henceforth, we assume $N$ is K\"{a}hler-Einstein.
The expression for $\triangle\kappa$, where the Ricci tensor of $N$ appears, and the Weitzenb\"{o}ck formula for $\Fw$, leading to an integral equation involving  the scalar curvature $R$, some trigonometric functions of 
the common K\"{a}hler angle, and the gradient 
of its cosine  (Proposition 4.2), are our tools to obtain the results of this paper.
In section 4 we prove our main results, namely:
\begin{Th}Let $F$ be a minimal immersion of a compact oriented manifold $M$,
into a K\"{a}hler-Einstein  manifold  $N$,
with  equal K\"{a}hler angles.\\[1mm]
(i)~~ If $n=2$ and $R\neq 0$, then  $F$ is either a complex 
or a Lagrangian submanifold.\\[1mm]
(ii)~~ If $n\geq 3$,  $R<0$, and $F$ has no complex points, 
 then $F$ is Lagrangian.\\[1mm]
(iii)~ If $n\geq 3$, $R=0$, and $F$ has no complex points,
 then the  common K\"{a}hler angle  must be 
constant. 
\end{Th}
The conclusions in $(i)$ and $(ii)$ give a generalization of Theorem 1.1  to higher dimensions
and equal K\"{a}hler angles.
The case $n=2$  is the most special, because, in this dimension,
 immersions with equal K\"{a}hler angles  have harmonic $\Fw$, as we will see in section 3. The case $n=3$  also has special properties. If the angle is constant we may allow $R>0$:
\begin{Th} Let $F$ be  minimal with constant equal K\"{a}hler angles,
$M$ compact, orientable, and  $R\neq 0$. Then,  $F$ is either a
complex  or a Lagrangian submanifold.
\end{Th}
\begin{Th} Let $F$ be  minimal with equal K\"{a}hler angles, and
M compact, orientable, with non-negative
 isotropic scalar curvature. If
$n=2$ or $3$,  then one of the 
following cases holds:\\[1mm]
$(i)$~~~ $M$ is a complex submanifold of $N$.\\
$(ii)$~~ $M$ is  a Lagrangian submanifold of $N$.\\
$(iii)$~ $R=0$ and $\cos\theta = constant \neq 0,1$, 
$J_{\omega}$ is a complex integrable structure, with
$(M,J_{\omega},g_{M})$ a K\"{a}hler manifold.\\[1mm]
For any  $n\geq 1$, any $R$,  and constant equal K\"{a}hler angle, 
$(i),(ii)$ or $(iii)$ hold
as well.
\end{Th}
This theorem can be applied, for instance, to flat minimal tori on 
Calabi-Yau manifolds, or to spheres or products of $S^{2}$ with $S^{2}$
or with flat tori minimaly immersed into K\"{a}hler-Einstein manifolds
with positive scalar curvature.
\section{The morphism $\Phi$}
\setcounter{Th}{0}
\setcounter{Pp}{0}
\setcounter{Cr} {0}
\setcounter{Lm} {0}
\setcounter{Def} {0}
\setcounter{equation} {0}
We consider the following morphism of vector bundles
\[
\begin{array}{cccc}
\Phi: &TM & \ra & NM \\
      & X & \ra & (JdF(X))^{\bot}\end{array}\\[-4mm]
\]
We easily verify that
\begin{equation}
\Phi(X)= JdF(X)- dF(F^{*}\omega (X)).
\end{equation}
Both $TM$ and $NM$ are real vector bundles of the same dimension $2n$. 
$F$ has no complex directions  iff $\Phi$ is an isomorphism.
In fact $\Phi(X)=0$, iff $JdF(X)=dF(Y)$ for some $Y$, i.e.,
$span\{X, Y=``JX"\}$ is a complex direction of $F$. Assume there are no
complex directions. Then, 
\begin{equation}
\hat{g}(X,Y)= g_M(X,Y)-g_M(F^{*}\omega (X), F^{*}\omega (Y))
\end{equation}
defines a Riemannian metric on $M$.
With this metric,  $\Phi:(TM,\hat{g})\ra (NM,g)$ 
 is an isomorphism of Riemannian vector bundles. Let us denote by
$\lnab{}$, $\hlnab{}$, $\lnabo{}$, and $\lnabp{}$, res\-pec\-ti\-vely, the
 Levi-Civita connection of $(M,g_{M})$, the Levi-Civita connection
of $(M,\hat{g})$, the usual connection of $NM$ induced by the
Levi-Civita connection of $N$, and the connection on $TM$ that makes
the isomorphism $\Phi$ parallel, namely $\lnabp{}= \Phi^{-1*}\lnabo{}$.
We will also denote by $\lnab{}$ the Levi-Civita connection of $N$
and the induced connection on $F^{-1}TN$, as well.
Thus, if $U$ is a smooth section of $NM\subset F^{-1}TN$, and $X,Y$ 
are smooth  vector fields on $M$, we have\\[-8mm]
\[\lnabo{X}U = (\lnab{X}U)^{\bot}~~~~~~
~~~~~~\Phi(\lnabp{X}Y)=\lnabo{X}(\Phi(Y)).\]
The connections $\lnab{}$ and $\hlnab{}$ have no torsion, because they 
are  Levi-Civita, but $\lnabp{}$ may have non-zero torsion $T'$.
Since both $\hlnab{}$ and $\lnabp{}$ are Riemannian connections 
of $TM$ for the same Riemannian metric $\hat{g}$, then $T'=0$
iff $\hlnab{}=\lnabp{}$ iff $\Phi$ is parallel. 
Note that, if $F$ is Lagrangian, then
$\Phi(X)=JdF(X)\in NM$,  $J(NM)= dF(TM)$, and
$\hat{g}=g_{M}$, $\hlnab{}=\lnab{}$. Therefore,
$\lnab{X}\Phi~(Y)= \la{(}J\lnab{X}dF(Y)\la{)}^\bot =0$, 
that is, $\Phi$ is parallel, and so $\lnabp{}=\lnab{}$, as well.
 In the next section (Corollary 3.2), we will see a converse of this.
We extend  $\Phi:TM^{c}\ra NM^{c}$  to the complexified
spaces by $\Co$-linearity.
\begin{Lm} If $\{X_{\al},Y_{\al}\}$ is a diagonalizing $g_{M}$-orthonormal  basis of
$F^{*}\omega$ at $p$, then at $p$, and for each $\al,\be$\\[-10mm]
\begin{eqnarray*} 
\Phi(T'(Z_{\al},Z_{\bbe})) & = &
i(\cos\theta_{\al}+\cos\theta_{\be})\lnab{Z_{\al}}dF(Z_{\bbe})\\
\Phi(T'(Z_{\al},Z_{\be})) & = &
i(\cos\theta_{\al}-\cos\theta_{\be})\lnab{Z_{\al}}dF(Z_{\be}).
\end{eqnarray*}
\end{Lm}
\em Proof. \em  \\[-10mm]
\begin{eqnarray*}
\Phi(\lnabp{X}Y) &= &\lnabo{X}(\Phi(Y))
=\la{(}\lnab{X}(\Phi(Y))\la{)}^{\bot}= \la{(}\lnab{X}(JdF(Y)-dF(F^{*}\omega(Y)))\la{)}^{\bot}\\
&=&\la{(} J\lnab{X}dF(Y) +JdF(\lnab{X}Y)-\lnab{X}dF(F^{*}\omega(Y))\la{)}^{\bot}.\\[-8mm]
\end{eqnarray*}
Therefore, using the symmetry of the $\lnab{} dF$ and the fact 
that $\lnab{}$ is torsionless,\\[-4mm]
\begin{equation}
\Phi(T'(X,Y)) = \Phi( \lnabp{X}Y-\lnabp{Y}X-[X,Y] ) =
-\lnab{X}dF(F^{*}\omega(Y)) + \lnab{Y}dF(F^{*}\omega(X)).
\end{equation}
The lemma follows now immediately.~~~~~~~~~~~~\qed\\[4mm]
For each $U\in NM_{p}$, let us denote by $A^{U}:T_{p}M\ra T_{p}M$ the
symmetric operator $g_{M}(A^{U}(X),Y)= g( \lnab{}dF(X,Y),U)$.
From Lemma 2.1 and (2.3) we have
\begin{Pp} If $F$ is an immersion without complex directions,  then:\\[1mm]
$(i)$~~~ $\Phi$ is parallel iff $F^{*}\omega$ anti-commutes with $A^{U}$,
$\forall U\in NM$.\\
$(ii)$~~ If $F$ has equal K\"{a}hler angles,  on $\Omega^{0}_{2n}$,
 $T'$ is of type $(1,1)$ with respect to $J_{\omega}$.\\
$(iii)$~ On  $\Omega^{0}_{2n}$, 
 $F$ is pluriminimal  iff $T'$ is of 
type $(2,0)+(0,2)$ with respect to $\Jw$.\\
$(iv)$~~ If $F$ is broadly-pluriminimal, then, for 
$p \in\Omega_{2k}^{0}$ with $k\geq 1$, $T'$ is of 
type $(2,0)+(0,2)$ with respect to  any $g_{M}$-orthogonal
complex structure  $\tilde{J}=\Jw\oplus J'$ on $T_{p}M$, 
 where $J'$ is any  $g_{M}$-orthogonal complex structure
of $\kw$. 
\end{Pp}
\em Remark 1. \em  ~~If we call $\omega_{NM}$ the restriction of the K\"{a}hler
form $\omega$ to the normal bundle $NM$, we see that, if $\{X_{\al},
Y_{\al}\}$ is a diagonalizing $g_{M}$-orthonormal basis of $\Fw$ at a point $p$, then
$\{ U_{\al}=\Phi(\frac{Y_{\al}}{\sin\theta_{\al}}), V_{\al}=
\Phi(\frac{X_{\al}}{\sin\theta_{\al}})\}$ is a diagonalizing 
$g$-orthonormal  basis of $\omega_{NM}$. Moreover, $NM$ has the 
same K\"{a}hler angles  as $F$. Let $J_{NM}$ denote the complex structure 
on $NM$ defined
by this basis, that is, the one that comes from the polar decomposition 
of $\omega_{NM}$. Then,
$\Phi\Jw=-J_{NM}\Phi$.\\[5mm]
We should also  remark the following:
\begin{Pp}  If $F$ is an immersion with parallel 2-form $F^{*}\omega$,
then the K\"{a}hler angles  are constant and, in particular,
$M=\Omega^{0}_{2k}$ for some $k$. In this case,
considering $TM$ with  the Levi-Civita connection $\lnab{}$,~
 ${\cal K}_{\omega}$ and
${\cal K}_{\omega}^{\bot}$ are parallel sub-vector bundles of $TM$,
 and $ J_{\omega}\in C^{\infty}({\cal K}_{\omega}^{\bot *}\otimes
{\cal K}^{\bot}_{\omega})$,
$\tilde{g},\, \hat{g}\in C^{\infty}(\bigodot^{2} T^{*}M)$
are parallel sections. Furthermore, $(X,Y,Z)\!\ra\!\gdf{Z}{X}{Y}$ is 
symmetric on $TM$,  and,  if $F$ has no complex directions,
 $\hlnab{}=\lnab{}$. Moreover, if $\cos\theta_{\al_{1}}>\ldots>
\cos\theta_{\al_{r}}$ are the distinct eigenvalues of $F^{*}\omega$,
the corresponding eigenspaces  $E_{\al_{t}}$ define a smooth integrable
distribution of $TM$ whose integral submanifolds are parallel submanifolds
of $M$. The integral submanifolds of $E_{\al_{r}}$ are isotropic in $N$ if
$\cos\theta_{\al_{r}}=0$, and the ones of $E_{\al_{1}}$ are complex
submanifolds of $N$ if $\cos\theta_{\al_{1}}=1$. The other ones
are K\"{a}hler manifolds with respect to $\Jw$, and $F$ restricted to each
one of them is an immersion of constant equal K\"{a}hler angles $\theta_{\al_{t}}$ with respect to $J$.
\end{Pp}
\em Proof. \em 
If $X,Y$ are smooth vector fields on $M$ and $Z\in T_{p}M$,
 an elementary computation gives
\begin{equation}
\lnab{Z}F^{*}\omega (X,Y)=-\gdf{Z}{X}{Y}+\gdf{Z}{Y}{X},
\end{equation}
which proves the symmetry of $(X,Y,Z)\ra\gdf{Z}{X}{Y}$.
 From (2.2) we see that $\hat{g}$ is parallel. 
Consequently, outside complex directions, $\lnab{}=\hlnab{}$.
If we parallel
transport a diagonalizing orthonormal basis $\{X_{\al},Y_{\al}\}$
 of $\Fw$ at $p_{0}$ along geodesics, on a
neighbourhood of $p_{0}$, since $\Fw$ is parallel we get
a diagonalizing orthonormal  frame on a whole neighbourhood with the property
$\lnab{}X_{\al}(p_{0})=\lnab{}Y_{\al}(p_{0})=0$. It also follows that
$\cos\theta_{\al}$ remains constant along geodesics, so it is constant,
and $\Jw(X_{\al})=Y_{\al}$ on a neighbourhood of $p_{0}$, with
$\lnab{}\Jw=0$ at $p_{0}$, and so $\Jw$ is parallel.
Si\-mi\-lar\-ly we see that $\tilde{g}$ is parallel. 
If we extend $F^*\omega$ to the complexified
tangent space $T^{c}_{p_{0}}M$, then 
$F^{*}\omega(Z_{\al})=i\cos\theta_{\al}
Z_{\al}$, and $F^{*}\omega(Z_{\bal})=-i\cos\theta_{\al}Z_{\bal}$.
Obviouly, the corresponding eigenspaces of $F^{*}\omega$,
are parallel sub-vector bundles of $T^{c}M$.
~~~~~~\qed\section{Immersions with equal K\"{a}hler angles }
\setcounter{Th}{0}
\setcounter{Pp}{0}
\setcounter{Cr} {0}
\setcounter{Lm} {0}
\setcounter{Def} {0}
\setcounter{equation} {0}
If $F$ has equal K\"{a}hler angles, then 
\[ F^{*}\omega = \cos\theta\,  J_{\omega}~~~~~\mbox{and}~~~~~\hat{g}=\sin^{2}\theta\,  g_{M},\]
with $\cos\theta$ a locally Lipschitz map on $M$,
smooth on the  open set where it does not vanish, and
$\Omega^{0}_{2k}=\emptyset$~ $\forall k\neq 0,n$. 
Note that $\sin^{2}\theta $ and $\cos^{2}\theta $ are smooth on all $M$.
The set
${\cal L}=\cos\theta^{-1}(\{ 0 \} )$ is the set of Lagrangian points,
for, at these points,  the tangent space of $M$ is a Lagrangian subspace
of the tangent space of $N$. Its subset of interior points is
$\Omega^{0}_{0}$. 
 Similarly, we say that 
${\cal C}=\cos\theta^{-1}(\{ 1 \} )$ is the set of complex points.
On the
open set $\Omega^{0}_{2n}= \cos\theta^{-1}(\R{}\sim \{ 0 \} )=M
\sim {\cal L}$, 
$J_{\omega}$ defines a smooth almost complex structure $g_{M}$-orthogonal. 
On the open set $\cos\theta^{-1}(\R{}\sim\{ 1 \} )=
M \sim {\cal C}$,  $\hat{g}$ is a
smooth metric conformally equivalent to $g_{M}$. Thus, if $n\geq 2$,
$~\hlnab{}=\lnab{}~$ iff $~\theta~$ is constant. 
Since the K\"{a}hler angles are equal, any smooth local
orthonormal  frame of the type $\{ X_{\al}, Y_{\al}=J_{\omega}X_{\al}\}$
diagonalizes $F^{*}\omega$ on the whole set where it is defined.
From  $\Fw=\cos\theta\Jw $, we get $\lnab{X}\Fw=
d\cos\theta(X)\Jw +\cos\theta \lnab{X}\Jw$, with  $\Jw$ orthogonal
to $\lnab{X}\Jw$ with respect to the Hilbert-Schmidt  inner product (because $\|\Jw\|^{2}=2n$ is constant).
Hence, considering $\Fw$  an operator on $TM$, on $\Omega^{0}_{2n}
\cup \Omega^{0}_{0}$
\begin{equation}
\|\lnab{}\Fw\|^{2} = 2n\|\nabla\cos\theta\|^{2}+
\cos^{2}\theta\|\lnab{}\Jw\|^{2}.
\end{equation}
We observe that $M\sim (\Omega^{0}_{2n}\cup \Omega^{0}_{0})$ is a set of Lagrangian points with no interior.
On $\Omega^{0}_{2n}$, we have then, 
$\lnab{}\Fw=0$ iff $\lnab{}\Jw=0$ and $\theta$ is constant.
Note that $\|\lnab{}\Fw\|^{2}$,
considering $\Fw$  an operator on $TM$, is twice the square norm when
considering $\Fw$  a 2-form.  From (2.3) we get, on $M\sim {\cal C}$,
\begin{equation}
\Phi(T'(X,Y))= 2\cos\theta (\lnab{}dF)^{(1,1)}(J_{\omega}X,Y).
\end{equation}
The right-hand side of (3.2) is defined to  be zero at a Lagrangian point.
Consequentely
\begin{Pp} If $F$ is an immersion with equal K\"{a}hler angles  and without complex points, then $T'=0$, that is, $\lnabp{}=\hlnab{}$ 
iff $\Phi$ is parallel iff $F$ is Lagrangian or pluriminimal.
In particular, if $F$ is minimal, $\Phi$ is parallel iff $F$ is 
broadly-pluriminimal.
\end{Pp}
 Let $Re(u+iv)=u$, for $u,v\in NM$.
\begin{Pp}
If $F$ is any immersion with equal K\"{a}hler angles, then, outside
complex and Lagrangian  points,\\[-5mm] \[
\Phi\la{(}\frac{1\!-\!n}{4}\nabla\log\sin^{2}\theta\la{)}=
 \frac{4\cos\theta}{\sin^{2}\theta}Re\!\left( i\sum_{\be,\mu}
\La{(}\gdf{\bmu}{\mu}{\be}\!-\!
\gdf{\bmu}{\be}{\mu}\La{)}\Phi(\bbe)\right),\\[-2mm]
\]
where $\nabla\log\sin^{2}\theta$ is the gradient with respect to $g_{M}$.
\end{Pp}
If $F$ is a complex submanifold on a open set, then $\Jw$ is the induced complex structure on $M$ and  $\lnab{}dF$ is of type $(2,0)$. Applying
Proposition 2.2 on $\Omega^{0}_{0}$, and Proposition 3.1 on open sets without
complex and Lagrangian points, and noting that $\{\Phi(\be), \Phi(\bbe)=\overline{\Phi(\be)}\}_
{1\leq \be\leq n}$ multiplied by $\frac{\sqrt{2}}{\sin\theta}$ 
constitutes an unitary basis of $NM^{c}$, we immediately conclude
\begin{Cr}
If $F$ is an immersion with equal K\"{a}hler angles, and $n\geq 2$, 
then  $\theta$ is  constant iff\\[-4mm]
\begin{equation}
\sum_{\mu}\gdf{\bmu}{\mu}{\be}
=\sum_{\mu}\gdf{\bmu}{\be}{\mu}~~~~~\forall\be.
\end{equation}
\end{Cr}
Note that  (3.3) is a sort of symmetry property, 
  and the first term is just $\frac{n}{2}g(H,JdF(\be))$, 
 where $H=\frac{1}{2n}trace_{g_{M}}\lnab{}dF=\frac{2}{n}
\sum_{\mu}\lnab{}dF(\bmu,\mu)$ is the mean curvature of $F$.
\begin{Th} If $n\geq 2$ and  $F$ is a pluriminimal immersion with equal K\"{a}hler angles  then $\cos\theta=$ constant. Moreover, if
it is not a complex submanifold,   then $\lnab{}=
\hlnab{}=\lnabp{}$, and $N$ must be Ricci-flat. In particular,
$\Phi$ defines a parallel homothetic isomorphism from
$(TM,g_{M})$ onto $(NM,g)$.
\end{Th} 
\em Proof. \em  On a neighbourhood of a non-complex point,
  from Proposition 3.1, $\hlnab{}=\lnabp{}$, and
from Corollary 3.1, $\cos\theta$ is constant. Then $\hlnab{}=\lnab{}$, 
as well. So if $F$ is not a complex submanifold, it has no complex
points anywhere. Finally, (1.7) for pluriminimal
immersions with $\kappa=constant$ gives $R=0$.~~~~~~~\qed\\[5mm]
The above theorem and Proposition 3.1 lead to:
\begin{Cr}  If $F$ is a minimal immersion with equal K\"{a}hler angles, 
 without complex
points, $n\geq 2$, and  $R\neq 0$, then $F$ is Lagrangian  iff $\Phi$
is parallel.
\end{Cr}
To prove Proposition 3.2 we will need to relate the three connections 
of $M$, $\lnab{}$, $\hlnab{}$ and $\lnabp{}$.  Let
$\{e_{1},\ldots,e_{2n}\}=
\{X_{\mu},Y_{\mu}=J_{\omega}X_{\mu}\}_{1\leq \mu \leq n}$ be
a local $g_{M}$-orthonormal  frame outside the Lagrangian and complex set,
and $\partial_{1},\ldots,\partial_{2n}$  a local frame of $M$ defined by
a coordinate chart. Set $g_{ij}=g_{M}(\partial_{i},\partial_{j})$,~
$\hat{g}_{ij}=\hat{g}(\partial_{i},\partial_{j})=\sin^{2}\theta g_{ij}$,
and $e_{s}=\sum_{i}\lambda_{si}\partial_{i}$.
The Christofel symbols are given by
$~~2\hat{\Gamma}^{k}_{ij} =\!\! \sum_{s}\hat{g}^{ks}(\partial_{i}
\hat{g}_{sj}+ \partial_{j}\hat{g}_{is}- \partial_{s}\hat{g}_{ij})
= \delta_{kj}\partial_{i}\log\sin^{2}\theta
+\delta_{ki}\partial_{j}\log\sin^{2}\theta-\!\sum_{s}g^{ks}g_{ij}
\partial_{s}\log\sin^{2}\theta +2\Gamma^{k}_{ij}.~$
Hence \\[-2mm]
\[
\hlnab{\partial_{i}}\partial_{j}-\lnab{\partial_{i}}\partial_{j}=
\sum_{k} ( \hat{\Gamma}^{k}_{ij}-\Gamma^{k}_{ij})\partial_{k}
= \frac{1}{2}\La{(}\partial_{i}(\log\sin^{2}\theta )\partial_{j}
+\partial_{j}(\log\sin^{2}\theta )\partial_{i}
-g_{ij}\nabla (\log\sin^{2}\theta) \La{)}\\[-1mm]
\]
Since  ~$\sum_{ij}g_{ij}\lambda_{si}\lambda_{sj}=1$,
~$\sum_{s}\hlnab{e_{s}}e_{s}
-\lnab{e_{s}}e_{s}=\sum_{sij}\lambda_{si}\lambda_{sj}(\hlnab{
\partial_{i}}\partial_{j} -\lnab{\partial_{i}}\partial_{j})=
(1-n)\nabla\log\sin^{2}\theta$.
Therefore,\\[-7mm]
\begin{eqnarray}
\lefteqn{\!\!\!\!\!\!\!\!\!\!\!\!\!\!\!\!\!\!\!\!\!\!\!\!\! 
\sum_{\mu}\!\hlnab{\bmu}\mu\!-\!\!\lnab{\bmu}\mu=
\frac{1}{4}\sum_{\mu}\!\!\La{(}
\!\hlnab{X_{\mu}}\!X_{\mu}\!+\!\!\hlnab{Y_{\mu}}Y_{\mu}\!
 -\!\!\lnab{X_{\mu}}X_{\mu} \!-\!\!\lnab{Y_{\mu}}Y_{\mu}\!\La{)}\!
-i\La{(}\!\!\hlnab{X_{\mu}}Y_{\mu}\!-\!\!\hlnab{Y_{\mu}}X_{\mu} \!-\!\!\lnab{X_{\mu}}Y_{\mu} 
\!+\!\lnab{Y_{\mu}}X_{\mu}\!\La{)}} \nonumber\\[-3mm]
&&\!\!\!\!\!\!\!\!\!\!\!\!\!\!\!\!\!\!\!\!\!\!\!\!\!
=\frac{1}{4}\sum_{s} (\hlnab{e_{s}}e_{s}
-\lnab{e_{s}}e_{s}) +
\frac{i}{4}\sum_{\mu}\La{(}[Y_{\mu},X_{\mu}]-
[Y_{\mu},X_{\mu}]\La{)} ~= ~~
\frac{(1-n)}{4}\nabla\log\sin^{2}\theta.
\end{eqnarray}
Set $ S'(X,Y)=\lnabp{X}Y-\hlnab{X}Y$.
Then $ S'(X,Y)- S'(Y,X)=T'(X,Y)$.
Similarly we get\\[-5mm]
\begin{equation}
 \sum_{\mu}\lnabp{\bmu}\mu
-\hlnab{\bmu}\mu =\frac{1}{4}trace_{g_{M}}S' -\frac{i}{4}\sum_{\mu}T'(X_{\mu},Y_{\mu}).\\[2mm]
\end{equation}
\begin{Lm} $\forall X\in T_{p}M$,~~  $\sum_{i}\hat{g}(S'(e_{i},e_{i}),X)
=-\sum_{i}\hat{g}(T'(e_{i},X),e_{i})$.\\[-5mm]
\end{Lm}
\em Proof. \em   We may assume that the local referencial $\partial_{i}$ is
$\hat{g}$-orthonormal  at a fixed poit $p_{0}$. On a neighbourhood of
$p_{0}$, we define  
 ${\Gamma'}_{ij}^{k}$ and ${S'}_{ij}^{k}$ as \\[-2mm]
\[\lnabp{\partial_{i}}\partial_{j}=
\sum_{k}{\Gamma'}_{ij}^{k}\partial_{k}
~~~~~~~~~~S'(\partial_{i},\partial_{j})=
\sum_{k}{S'}_{ij}^{k}\partial_{k}=\sum_{k}
({\Gamma'}_{ij}^{k}-\hat{\Gamma}_{ij}^{k})\partial_{k}.\\[-2mm]\]
Then $~{T'_{ij}}^{k}={\Gamma'}_{ij}^{k}-{\Gamma'}_{ji}^{k}~$, and,
 at $p_{0}$,~ ${\Gamma'}_{ij}^{k}=\hat{g}(\lnabp{\partial_{i}}\partial_{j},\partial_{k}),
~~~{S'}_{ij}^{k}=\hat{g}(S'(\partial_{i},\partial_{j}),\partial_{k})=
{\Gamma'}_{ij}^{k}-\hat{\Gamma}_{ij}^{k}$.
 $\lnabp{}$ is a Riemannian connection with respect to $\hat{g}$.
Then
\[
\partial_{i}\hat{g}_{jk}(p_{0})=
\hat{g}(\lnabp{\partial_{i}}\partial_{j},
\partial_{k})+\hat{g}(\partial_{j},\lnabp{\partial_{i}}\partial_{k})=
{\Gamma'}_{ij}^{k}+{\Gamma'}_{ik}^{j}\\[-3mm]
\]
Hence, at $p_{0}$
\begin{eqnarray*}
2\hat{\Gamma}_{ij}^{k}\!&=&\!\!\sum_{s}\hat{g}^{ks}(\partial_{i}
\hat{g}_{sj}+ \partial_{j}\hat{g}_{is}- \partial_{s}\hat{g}_{ij})
={\Gamma'}_{ik}^{j}+{\Gamma'}_{ij}^{k}+
{\Gamma'}_{ji}^{k}+{\Gamma'}_{jk}^{i}-{\Gamma'}_{ki}^{j}
-{\Gamma'}_{kj}^{i}\\[-2mm]
\!&=&\!\!({\Gamma'}_{ij}^{k}+{\Gamma'}_{ji}^{k})\!+\!
({\Gamma'}_{ik}^{j}-{\Gamma'}_{ki}^{j})\!
+\!({\Gamma'}_{jk}^{i}-{\Gamma'}_{kj}^{i})=
({\Gamma'}_{ij}^{k}+{\Gamma'}_{ji}^{k})+
{T'}_{ik}^{j}+{T'}_{jk}^{i}
\end{eqnarray*}
But $~{\Gamma'}_{ij}^{k}+{\Gamma'}_{ji}^{k}=
2{\Gamma'}_{ij}^{k}+({\Gamma'}_{ji}^{k}-{\Gamma'}_{ij}^{k})
=2{\Gamma'}_{ij}^{k}+{T'}_{ji}^{k}~$.
Thus
\[
{S'}_{ij}^{k}={\Gamma'}_{ij}^{k}-\hat{\Gamma}_{ij}^{k}=
\frac{1}{2}({T'}_{ij}^{k}-{T'}_{ik}^{j}+{T'}_{kj}^{i}). 
\]
That is, at $p_{0}$,~
$\hat{g}(S'(\partial_{i},\partial_{j}),\partial_{k})=
\frac{1}{2}\La{(}\hat{g}(T'(\partial_{i},\partial_{j}),\partial_{k})
-\hat{g}(T'(\partial_{i},\partial_{k}),\partial_{j})
+\hat{g}(T'(\partial_{k},\partial_{j}),\partial_{i})\La{)}.$
We may assume that, at $p_{0}$, 
$\partial_{i}(p_{0})=\frac{e_{i}}{\sin\theta}$, leading to the Lemma.~~~~~~~~~~~\qed\\[7mm]
\em Proof of  Proposition 3.2. \em  ~~~Following  the proof of Lemma 2.1,
$~~\Phi(\lnabp{X}\mu -\lnab{X}\mu)=$ \newline $=
 ((J-i\cos\theta)\lnab{X}dF(\mu))^{\bot}~$.
Hence, from (3.4),\\[-4mm]
\[
\Phi(\frac{(1\!-\!n)}{4}\nabla\log\sin^{2}\theta)=
\Phi(\sum_{\mu}\hlnab{\bmu}\mu -\! \lnab{\bmu}\mu)
=\LA{(}\!(J-i\cos\theta)\frac{nH}{2}\!\LA{)}^{\bot}\!\!-\!
\sum_{\mu}\Phi(\lnabp{\bmu}\mu -\! \hlnab{\bmu}\mu).
\]
But, from (3.5),~ $\sum_{\mu}\Phi(\lnabp{\bmu}\mu - \hlnab{\bmu}\mu)=
\frac{1}{4}\Phi(trace_{g_{M}}S')-\frac{i}{4}\Phi(\sum_{\mu}
T'(X_{\mu},Y_{\mu}))$.~
The skew-symmetry of $T'$ and (3.2) implies that
\[ \Phi(\sum_{\mu}T'(X_{\mu},Y_{\mu}))= -2i\sum_{\mu}\Phi(T'(\mu,\bmu))
=4\cos\theta\lnab{\mu}dF(\bmu)= 2n\cos\theta H.\]\\[-5mm]
Thus, ~$\sum_{\mu}\Phi(\lnabp{\bmu}\mu-\hlnab{\bmu}\mu)=\frac{1}{4}
\Phi(trace_{g_{M}}S') -\frac{ni}{2}\cos\theta H.$~
Therefore,
\begin{equation}
\Phi(\frac{(1-n)}{4}\nabla\log\sin^{2}\theta)=
\frac{1}{4}\LA{(}2n(JH)^{\bot}-\Phi(Trace_{g_{M}} S') \LA{)}.
\end{equation}
Using Lemma 3.1,  (3.2),  and 
$\Phi(\mu)=JdF(\mu)-i\cos\theta dF(\mu)$,
we have
\begin{eqnarray*}
&&\!\!\!\!\!\!\!\Phi(Trace_{g_{M}} S')
=\sum_{j,k}\hat{g}\La{(}S'(e_{j},e_{j}), \frac{e_{k}}{\sin\theta}\La{)}
\Phi(\frac{e_{k}}{\sin\theta})
=\sum_{j,k}-\hat{g}\La{(}T'(e_{j},\frac{e_{k}}{\sin\theta}),e_{j}\La{)}
\Phi(\frac{e_{k}}{\sin\theta})\\[-2mm]
&&\!\!\!\!\!\!\!
=\frac{-4}{\sin^{2}\theta}\!\!\sum_{\mu,\be}\!
\LA{(}\!\La{(}\hat{g}(T'(\mu,\be),\bmu)
\!+\!\hat{g}(T'(\bmu,\be),\mu)\!\La{)}\Phi(\bbe)
+\La{(}\hat{g}(T'(\mu,\bbe),\bmu)
\!+\!\hat{g}(T'(\bmu,\bbe),\mu)\!\La{)}\Phi(\be)\!\LA{)}\\[-2mm]
&&\!\!\!\!\!\!\!=-\frac{4}{\sin^{2}\theta}\sum_{\mu,\be}\LA{(}
g\La{(}\Phi\la{(}T'(\bmu,\be)\la{)},\Phi(\mu)\La{)}\Phi(\bbe)
+g(\Phi\la{(}T'(\mu,\bbe)\la{)},\Phi(\bmu)\La{)}\Phi(\be)~\LA{)}\\[-2mm]
&&\!\!\!\!\!\!\!=\frac{8i\cos\theta}{\sin^{2}\theta}
\sum_{\mu,\be}\LA{(}\gdf{\bmu}{\be}{\mu}\Phi(\bbe)
-\gdf{\mu}{\bbe}{\bmu}\Phi(\be)\LA{)}. 
\end{eqnarray*}
Writing $(JH)^{\bot}$ in terms of $\Phi(\be)$ and $\Phi(\bbe)$,
\begin{eqnarray*}
2n(JH)^{\bot}&=& \sum_{\be}\frac{4n}{\sin^{2}\theta}
\La{(}g(JH,\Phi(\be))\Phi(\bbe)+g(JH,\Phi(\bbe))\Phi(\be)\La{)}\\[-3mm]
&=&\sum_{\be,\mu}\frac{8i\cos\theta}{\sin^{2}\theta}\LA{(}
\gdf{\bmu}{\mu}{\be}\Phi(\bbe) - \gdf{\bmu}{\mu}{\bbe}\Phi(\be)\LA{)},
\end{eqnarray*}
and substituing these equations into (3.6),
 we prove Proposition 3.2.~~~~~~\qed

\subsection{The Weitzenb\"{o}ck formula for $\Fw$}
For simplicity let us use the  notation
\[ \gf{X}{Y}{Z}=\gdf{X}{Y}{Z}.\]
We also observe that, from \\[-5mm]
\begin{equation}
\forall \mu~~~~  \frac{i}{2}\cos\theta = F^{*}\omega(\mu,\bmu),
\end{equation}
valid on an open set, and from (2.4), we obtain $\forall \mu$
\begin{eqnarray}
 \frac{i}{2}d\cos\theta (X) &= &d(\Fw (\mu,\bmu) )(X)=
\lnab{X}\Fw (\mu,\bmu)+\Fw(\lnab{X}\mu,\bmu)+\Fw(\mu,\lnab{X}\bmu)\non\\[-2mm]
&=&-\gf{X}{\mu}{\bmu}+\gf{X}{\bmu}{\mu}+
2(\langle\lnab{X}\mu,\bmu\rangle+
\langle\lnab{X}\bmu,\mu\rangle)\Fw(\mu,\bmu)\non\\[-1mm]
&=&-\gf{X}{\mu}{\bmu}+\gf{X}{\bmu}{\mu}\mbox{~~~~~~(no~sumation~over~}\mu).
\end{eqnarray}
Then (3.3) is equivalent to  $\gdf{X}{\mu}{\bmu}=
\gdf{X}{\bmu}{\mu},~\forall\mu $ (or some $\mu$). 
 From
$\Jw Z_{\al}=iZ_{\al}$, $\Jw Z_{\bal}=-iZ_{\bal}$ and the fact that
$\Jw$ is $g_{M}$-orthogonal, we get, on $\Omega^{0}_{2n}$, $\forall
\al,\be$, and $\forall v\in TM$\\[-4mm]
\begin{equation}
\langle \lnab{v} \Jw (\al), \be\rangle = 
2i\langle \lnab{v} \al ,\be \rangle,~~~~~~~~~~~~
\langle \lnab{v} \Jw (\al), \bbe\rangle = 0.
\end{equation}
Recall that, if $\xi$ is a $r\!+\!1$-form on $M$, $r\geq 0$,  with
values on a vector bundle $E$ over $M$ with a connection $\lnab{}^{E}$,
then $\delta \xi$, the divergence of $\xi$,  is the $r$-form on $M$
with values on $E$ given by 
\[\delta\xi(u_{1},\ldots,u_{r})=-\sum_{s}\lnab{e_{s}}^{E}\xi(e_{s},
u_{1},\ldots,u_{r}),\]
where $e_{1},\ldots,e_{m}$ is an orthonormal  basis of $T_{p}M$,
$u_{i}\in T_{p}M$, and $\lnab{}^{E}\xi$ is the covariant derivative of 
$\xi$ on $\bigwedge^{r+1}T^{*}M\otimes E$.
 Thus,  $\delta$ is the formal adjoint of
$d$ on forms (cf. [E-L]). Note that $\delta \Fw (X)=\langle \delta \Fw, X\rangle, ~\forall X\in T_{p}M$,  considering on the  left-hand side $\Fw$ a (closed) 2-form
on $M$ and on the right-hand side   an endomorphism of $TM$.
\begin{Pp} Let $F$ be an immersion with equal K\"{a}hler angles   and 
$\nabla \cos\theta$  denote the gradient with respect to $g_{M}$. 
On $\Omega^{0}_{2n}$, and considering $F^{*}\omega$  an endomorphism of $TM$.
\\[-4mm]
\[ \delta F^{*}\omega= (n-2)J_{\omega}(\nabla\cos\theta),~~~~~~~~~~
\cos\theta(\delta J_{\omega})=(n-1)J_{\omega}(\nabla\cos\theta).\]\\[-8mm]
Thus,\\
$(i)$~~ For $n=1$, ~$\delta  J_{\omega}=0 $ (obviously!),
and $~\delta F^{*}\omega =0$ iff $\theta$ is constant.\\
$(ii)$~ For $n=2$,  ~$\delta F^{*}\omega=0$. Moreover,
$\delta J_{\omega}=0$ iff $\theta$ is constant.\\
$(iii)$~ For $n\neq 1,2$,~
 $\delta F^{*}\omega=0$ iff $\delta J_{\omega}=0$
 iff $\theta$ is constant.
\end{Pp}
\em Proof. \em   Considering $\Fw$  a 2-form on $M$, using the symmetry
of $\lnab{}dF$ and (2.4), if $X\!\in\! T_{p}M$,
\begin{eqnarray*}
\delta (\Fw)(X)=&&\!\!\!\!\!\!\!\!\!\sum_{\mu}\!\!
-2\lnab{\mu}\Fw (\bmu,X)-2\lnab{\bmu}\Fw (\mu,X)
= \sum_{\mu}2\gf{\mu}{\bmu}{X}\!-\!2\gf{\mu}{X}{\bmu}
\!+\!2\gf{\bmu}{\mu}{X}\!-\!2\!\gf{\bmu}{X}{\mu}\\[-3mm]
=&&\!\!\!\!\!\! \!\!\!2\sum_{\mu}(
-\gf{X}{\mu}{\bmu}+\gf{X}{\bmu}{\mu})
-4\sum_{\mu}(\gf{\bmu}{X}{\mu}-\gf{\bmu}{\mu}{X}).
\end{eqnarray*}
From (3.8), $\frac{ni}{2}d\cos\theta (X)=\sum_{\mu} -\gf{X}{\mu}{\bmu}+
\gf{X}{\bmu}{\mu}$. Therefore,
\begin{equation}
\delta (\Fw) (X) =
nid\cos\theta (X)-4\sum_{\mu}\lnab{\bmu}\Fw (\mu,X).
\end{equation}
Since  $\Fw$ is of type $(1,1)$ with respect to $\Jw$, and 
$\forall Z\in T_{p}^{c}M$, $\forall \mu,\be$, 
$\langle \lnab{Z}{\be},\mu\rangle = -\langle \be ,\lnab{Z}
\mu \rangle$, we get using (3.9)
\begin{eqnarray}
\lnab{Z}\Fw (\mu,\be)&=&d(\Fw(\mu,\be))(Z)-
\Fw(\lnab{Z}\mu,\be)-\Fw(\mu,\lnab{Z}\be)\non\\
&=&2i\cos\theta \langle \lnab{Z}\mu,\be \rangle =
\cos\theta \langle \lnab{Z}\Jw (\mu),\be \rangle.
\end{eqnarray}
 Note that, since $\Jw^{2}=-Id$, 
 $\lnab{X}\Jw (\Jw Y) = -\Jw ( \lnab{X}\Jw (Y))$, 
$\forall X,Y\in T_{p}M$. So
\begin{eqnarray*}
4\sum_{\mu}\lnab{\bmu}\Jw (\mu)&=&\sum_{\mu}
\lnab{X_{\mu}}\Jw (X_{\mu})+\lnab{Y_{\mu}}\Jw (Y_{\mu})
+i\lnab{Y_{\mu}}\Jw (X_{\mu})-i\lnab{X_{\mu}}\Jw (Y_{\mu})\\[-3mm]
&=&-\delta \Jw + i\sum_{\mu}(-\lnab{X_{\mu}}\Jw (\Jw X_{\mu}) -
\lnab{Y_{\mu}}\Jw (\Jw Y_{\mu}) )
=-\La{(}\delta \Jw +i\Jw(\delta\Jw)\La{)}.
\end{eqnarray*}
Hence, from (3.11), and since $\Jw$ is $g_{M}$-orthogonal, $\forall \be$\\[-4mm]
\[\sum_{\mu}\lnab{\bmu}\Fw (\mu,\be)=-\frac{\cos\theta}{4}\langle 
\delta\Jw +i\Jw(\delta\Jw),\be\rangle=-\frac{\cos\theta}{2}\langle 
\delta\Jw ,\be\rangle.\]
Moreover, $~id\cos\theta (\be)=d\cos\theta(\Jw\be)=\langle
\nabla\cos\theta,\Jw\be\rangle=-\langle \Jw(\nabla\cos\theta),\be\rangle.$~
From (3.10), 
$\delta\Fw (\be)=\langle -n\Jw(\nabla\cos\theta)+2\cos\theta\,
\delta \Jw~, ~\be\rangle$.
Thus,  if we consider $\Fw$  an endomorphism of $TM$, and
since $\langle,\rangle$, $\Jw$, and $\Fw$ are real operators,
\begin{equation}
\delta\Fw = -n\Jw(\nabla\cos\theta) +2\cos\theta\,\delta\Jw.
\end{equation}
On the other hand,  $\Fw=\cos\theta\Jw$. 
Then, an elementary computation gives
\begin{equation}
\delta\Fw = -\Jw(\nabla\cos\theta) +\cos\theta\,\delta\Jw.
\end{equation}
Comparing (3.12) with (3.13) we get the Proposition.~~~~~~~~~~~\qed
\\[4mm]
\em Remark 2. \em One may check the equation in Proposition 3.2 by using the  equalities given in the above Proposition and its proof.\\[5mm]

If we apply the Weitzenb\"{o}ck formula to the 2-form $\Fw$, for an immersion
$F:M\ra N$ we get (see e.g [E-L] (1.32))
\begin{equation}
\frac{1}{2}\triangle \|\Fw\|^{2} =-\langle \triangle \Fw,\Fw\rangle
+\|\lnab{}\Fw\|^{2}+\langle S\Fw, \Fw\rangle,
\end{equation}
where $\langle,\rangle$ denotes the Hilbert-Schmidt inner product for 2-forms,
and $S$ is the Ricci operator of $\bigwedge^{2}T^{*}M$. We note that we use
the the  sign convention $\triangle \phi = +Trace_{g_{M}}Hess\,\phi$,
for $\phi$ a smooth real map on $M$. This sign is opposite to the one of [E-L], but here we use the same sign as in [E-L] for the Laplacian of 
forms $\triangle = d\delta +\delta d$. If
$\overline{R}$ denotes the curvature tensor of $\bigwedge^{2}T^{*}M$,
and $X,Y,u,v\in T_{p}M$, $\xi\in \bigwedge^{2}T^{*}_{p}M$, then
\begin{eqnarray*}
&&\overline{R}(X,Y)\xi ~(u,v) = -\xi(R^{M}(X,Y)u\,,\, v)
-\xi(u\,,\,R^{M}(X,Y)v),\\
&&S\Fw(X,Y)=\sum_{1\leq i\leq 2n}-\overline{R}(e_{i},X)\Fw~(e_{i},Y) +
\overline{R}(e_{i},Y)\Fw~(e_{i},X),
\end{eqnarray*}
Where $R^{M}$ denotes the curvature tensor of $M$. In general, we use the
following sign convention for curvature tensors: $R^{M}(X,Y)Z=
-\lnab{X}\lnab{Y}Z+\lnab{Y}\lnab{X}Z+\lnab{[X,Y]}Z$. Then, $R^{M}(X,Y,Z,W)
=g_{M}(R^{M}(X,Y)Z,W)$.
It is straightforward to prove
\begin{Lm} If $\{X_{\al},Y_{\al}\}$ is a diagonalizing orthonormal basis
of $\Fw$ at $p$,
\begin{eqnarray*}
\langle S\Fw,\Fw\rangle\!\!
&=&\!\! \sum_{\mu}4\cos^{2}\theta_{\mu}Ricci^{M}(\mu,\bmu)+\sum_{\mu,\rho}
8\cos\theta_{\mu}\cos\theta_{\rho}R^{M}(\rho,\brho,\mu,\bmu)\\[-2mm]
 &=&\!\!\sum_{\mu,\rho}4
(\cos\theta_{\mu}\!+\!\cos\theta_{\rho})^{2}R^{M}(\rho,\mu,\brho,\bmu)
+4(\cos\theta_{\mu}\!-\!\cos\theta_{\rho})^{2}R^{M}(\brho,\mu,\rho,\bmu).
\end{eqnarray*}
In particular,
 if $F$ has equal K\"{a}hler angles  at $p$,
then, at $p$,\\[-5mm]
 \[\langle S\Fw,\Fw\rangle =16\cos^{2}\theta\sum_{\rho,\mu}
R^{M}(\rho,\mu,\brho,\bmu).\\[-3mm]\]
Moreover, if $(M,\Jw,g_{M})$ is
K\"{a}hler in a neighbourhood of $p$, then $\langle S\Fw,\Fw\rangle=0$.
\end{Lm}
For example, if $M$ has constant sectional curvature  $K$, 
$~\langle S\Fw,\Fw\rangle = 4(n-1)K\|\Fw\|^{2}$.
 If $(M,\Jw,g_{M})$ is a K\"{a}hler manifold
of constant holomorphic sectional  curvature $K$ then
$\langle S\Fw,\Fw\rangle $ $=
4K\LA{(}n\sum_{\mu}\cos^{2}\theta_{\mu}-\la{(}\sum_{\mu}\cos\theta_{\mu}
\la{)}^{2}\LA{)}$ has constant sign, with equality to 
zero iff $K=0$ or $F$ has equal K\"{a}hler angles.
If $\Fw$ is parallel, from 
 (3.14), we obtain that  $~\langle S\Fw,\Fw\rangle=0$. In the latter
case, if $n\geq 2$ and $M$ has constant sectional curvature, then, 
either $F$ is Lagrangian, or $K=0$.\\[2mm]

We recall the concept of non-negative \em isotropic
sectional curvature, \em for $M$ with dimension $\geq 4$,
defined  by Micallef and Moore in [Mi-Mo].  Let
\[ K_{\footnotesize isot}(\sigma)=\frac{ R^{M}(z, w,\bar{z}, \bar{w})}
{||z\wedge w||^{2}},\]
where $\sigma=span_{\Co}\{z,w\}$ is a totally isotropic complex
 two-plane in 
$T^{c}M$, that is, $u\in \sigma \Rightarrow  g_{M}(u,u)=0$, and where 
$R^{M}(x,y,u,v)$ is  extend
to the complexified tangent space by $\Co$-multilinearity. The curvature
of $M$ is said to be non-negative (resp.~positive) on totally
isotropic two-planes at $p$,  if $K(\sigma)\geq 0$ (resp.  $>0$) 
whenever
$\sigma \subset T_{p}^{c}M$ is a totally isotropic two-plane over $p$.
If $M$ is compact, simply connected and has positive isotropic sectional curvature everywhere, then $M$ is homeomorphic to a sphere ([Mi-Mo]). 
If $n\geq 1$, $T^{2n}$ is the flat torus, and $S^{2}$ is the euclidean
sphere of $\R{3}$, then $S^{2}\times T^{2n}$, $S^{2}\times S^{2}$,
$S^{2}\times S^{2}\times T^{2n}$ have isotropic sectional
curvature $\geq 0$ but not $>0$.
If $\{X_{\al},Y_{\al}\}$ is any orthonormal  basis of $T_{p}M$, and $``\mu"$ 
denotes $Z_{\mu}$ as in (1.3),  the expression
\begin{equation}
S_{isot}(\{Z_\al\}_{1\leq \al \leq n})=
\sum_{\rho\neq \mu} K_{isot}(span_{\Co}\{\rho,\mu\})=
4\sum_{\rho,\mu}R^{M}(\rho,\mu,\brho,\bmu)
\end{equation}
is a hermitian trace of the curvature of $M$ restricted to the maximal totally
isotropic subspace $span_{\Co}\{Z_{1},\ldots,Z_{n}\}$ of $T^{c}M$. To require it to be $\geq 0$, for all maximal totally isotropic subspaces - 
 and we will say that $M$ has non-negative \em isotropic scalar curvature \em
- seems to be  strictly weaker than to have non-negative isotropic
sectional curvature.   We also note that, any other metric conformaly equivalent  to the flat metric $g_{0}$ on the 2n-torus with non-negative isotropic scalar curvature  is homothetically equivalent to $g_{0}$, hence
 flat. In fact, in general, if $\hat{g}=e^{\phi}g_{M}$ is a conformaly equivalent metric on $M$, then, for each $g_{M}$-orthonormal basis 
$\{X_{\al},Y_{\al}\}$, 
$ \hat{S}_{isot}(\{\hat{Z}_{\al}\})=
e^{-\phi}S_{isot}(\{Z_\al\})- (n-1)e^{-2\phi}(2\triangle \phi + (n-1)\|\nabla\phi\|^{2})$, 
where $\hat{Z_{\al}}$ are defined by the $\hat{g}$-orthonormal basis
$\{e^{-\frac{\phi}{2}}X_{\al}, e^{-\frac{\phi}{2}}Y_{\al}\}$. To require 
$2\triangle \phi + (n-1)\|\nabla\phi\|^{2}\leq 0$,
implies, in case of $M$ compact, $\phi$ constant. We  observe that,
if $dim_{\R{}}M\geq 6$, then $S_{isot}\equiv 0$ does not imply
$M$ to be flat, but $K_{isot}\equiv 0$ implies so. We also note that,
if $dim_{\R{}}(T_{p}M)=4$, the set of curvature operators at $p$
with zero isotropic sectional curvature, is a vector space of dimension 9.\\[7mm]
Recall that, for an immersion with equal K\"{a}hler angles, 
  $\Fw$ is parallel iff
$\theta$ is constant and if $\cos\theta\neq 0$,
 $(M,\Jw, g_{M})$ is a K\"{a}hler manifold.
We are going to see that an extra condition on the scalar isotropic curvature of $M$ may imply
itself that the K\"{a}hler angle   is constant and/or $\lnab{}\Jw=0$.  From
Proposition 3.3, for any $n\geq 1$, on $\Omega^{0}_{2n}\cup \Omega^{0}_{0}$
\begin{equation}
\|\delta\Fw\|^{2}=(n-2)^{2}\|\nabla\cos\theta\|^{2}.
\end{equation}
In particular, if $n\neq 2$, $\|\nabla\cos\theta\|^{2}$ can be extended
as a smooth map on all $M$ (recall that $\Omega^{0}_{2n}\cup \Omega^{0}_{0}$
is dense on $M$), and from (3.1) we  get that $\cos^2\theta\|\lnab{}\Jw\|^{2}$
is also smooth.
Observe that $\|\delta\Fw\|^{2}$ has the same value
considering $\delta\Fw$  a vector or  a 1-form, but 
considering $\Fw$  a 2-form (as in (3.14))
 $\|\lnab{}\Fw\|^{2}$ is half of the square norm when
considering $\Fw$  an operator of $TM$ (as in (3.1)). For $n=2$,
$\Fw$ is co-closed, and so it is a harmonic 2-form. In fact, since $F$
has equal K\"{a}hler angles,  $\Fw=\cos\theta (X_{*}^{1}\wedge
Y_{*}^{1} +X_{*}^{2}\wedge Y_{*}^{2})$, and so $\ast \Fw =\pm \Fw$, where
$\ast$ is the Hodge star-operator of $(M,g)$, and
the $\pm$ sign depends on the orientation of the diagonalizing basis. 
In particular, $\Fw$ is co-closed. For $n\geq 3$, $\Fw$ is harmonic
iff $\theta$ is constant.

Integrating 
(3.14) on $M$, using (3.16) and (3.1), and the fact that $\int_{M}\langle\triangle\Fw,\Fw\rangle Vol_{M}
=\int_{M}\|\delta\Fw\|^{2}Vol_{M}$, we have\\
\begin{equation}
0=\int_{M}\LA{(}(n-(n-2)^2)\|\nabla\cos\theta\|^{2} +\frac{1}{2}
\cos^2\theta\|\lnab{}\Jw\|^{2}\LA{)}Vol_{M} +\int_{M}\langle S\Fw,
\Fw\rangle Vol_{M}.\\[3mm]
\end{equation}
The first integrand is smooth on $M$, for all $n$ ( for $n\!=\!2$
it gives half of (3.1)). The factor $n\!-\!(n\!-\!2)^2$ is respectively, $>\!0$, $=\!0$, $<\!0$, according $n=2$ or $3$, $n=4$, and
 $n\geq 5$. If $M$ has non-negative isotropic scalar curvature,  $\langle S\Fw,\Fw\rangle \geq0 $, by Lemma 3.2. 
We conclude:
\begin{Pp} Let $F$ be a  non-Lagrangian immersion with equal K\"{a}hler angles of a compact
 orientable $M$ with non-negative isotropic scalar curvature into a 
K\"{a}hler manifold $N$.
If $n=2$ or $3$, then $\theta$ is constant and $(M,\Jw,g_{M})$ is a K\"{a}hler manifold.  If $n=4$, $(\Omega^{0}_{2n},\Jw,g_{M})$ is a K\"{a}hler manifold
 (but $\theta$ does not need to be constant).
 For any $n\geq 1$ and $\theta$ constant, 
$\Fw$ is parallel, i.e.,  $(M,\Jw,g_{M})$ is a K\"{a}hler manifold.
\end{Pp}
\section{Minimal immersions with equal K\"{a}hler angles}
\setcounter{Th}{0}
\setcounter{Pp}{0}
\setcounter{Cr} {0}
\setcounter{Lm} {0}
\setcounter{Def} {0}
\setcounter{equation} {0}
Let us assume that $F:M\ra N$ is minimal with equal K\"{a}hler angles. On a open set of $M\sim {\cal L}$ where  a  orthonormal  frame
$\{X_{\al}, Y_{\al}=\Jw (X_{\al})\}$ is defined, we have from (3.11)
and (2.4),
for any $p$, $Z\in T_{p}M$ and $\mu,\ga$,\\[-3mm]
\begin{equation}
2\cos\theta \langle \lnab{Z}\mu,\ga\rangle
=-i\lnab{Z}\Fw (\mu,\ga) =ig_{Z}\mu\ga-ig_{Z}\ga\mu.
\end{equation}
Note that  $\Fw (\lnab{Z}\mu, \bga) $ $=
i\cos\theta \langle \lnab{Z}\mu,\bga\rangle$ 
$=-i\cos\theta \langle \mu,\lnab{Z}\bga\rangle$ $=-\Fw(\mu,\lnab{Z}\bga)$.
Hence, if $\mu\neq \ga$, 
$\lnab{Z}\Fw (\mu,\bga)=d(\Fw(\mu,\bga))(Z)=0$. Thus \\[-3mm]
\begin{equation}
g_{Z}\mu\bga=g_{Z}\bga\mu,~~~~~\forall\mu\neq\ga\\[-4mm]
\end{equation} \\[-3mm]
From (3.8), for each $\mu$, \\[-4mm]
\begin{equation}
-\frac{i}{2}d\cos\theta(Z)=-\lnab{Z}\Fw (\mu,\bmu)
= g_{Z}\mu\bmu -g_{Z}\bmu\mu ~~(\mbox{no~sumation~over~}\mu)
\end{equation} \\[-6mm] 
From (1.6), on $M\sim ({\cal L}\cup {\cal C})$ 
\begin{eqnarray}
\triangle \kappa
&=&4i\sum_{\be} Ricci^{N}(JdF({\be}),dF({\bbe}))\non\\[-2mm]
&&+\frac{32}{\sin^{2}\theta}\sum_{\be,\mu}\!Im \La{(}
R^{N}(dF({\be}),dF({\mu}),dF({\bbe}),
JdF({\bmu})\!+\!i\cos\theta dF(\bmu))\La{)}\\[-2mm]
&&-\frac{128\cos\theta}{\sin^{4}\theta}
\sum_{\be,\mu,\rho}Re\la{(}g_{\be}{\mu}{\brho}\,g_{\bbe}{\rho}{\bmu}\la{)}
\\[-2mm]
&&+\frac{64\cos\theta}
{\sin^{2}\theta}\sum_{\be,\mu,\rho}\LA{(}
|\langle\lnab{{\be}}{\mu},{\rho}\rangle|^{2} +
|\langle\lnab{{\bbe}}{\mu},{\rho}\rangle|^{2}\LA{)},
\end{eqnarray}
where now $\kappa= n\log \left( \frac{1+\cos\theta}{1-\cos\theta}\right)$.
Since $R(X,Y,Z,JW)$ is skew-symmetric on $(X,Y)$ and symmetric on $(Z,W)$,
$\sum_{\mu,\be}R^{N}(dF(\be),dF(\mu),dF(\bbe),JdF(\bmu))=0$.
Then, from the Gauss equation and  minimality of $F$, 
\begin{eqnarray*}
(4.4)&=&\sum_{\be,\mu}\frac{32}{\sin^{2}\theta}Im 
\La{(}i\cos\theta R^{N}(dF({\be}),dF({\mu}),dF({\bbe}),
dF(\bmu))\La{)}\\[-2mm]
&=&\frac{32\cos\theta}{\sin^{2}\theta}\sum_{\be,\mu}
R^{M}(\be,\mu,\bbe,\bmu)+ 
g(\lnab{}dF(\be,\bmu),\lnab{}dF(\mu,\bbe)).
\end{eqnarray*}
Using the unitary basis $\{\frac{\sqrt{2}}{\sin\theta}\Phi(\rho),
\frac{\sqrt{2}}{\sin\theta}\Phi(\brho)\}$ of the normal bundle,
\begin{eqnarray}
\lefteqn{\frac{32\cos\theta}{\sin^{2}\theta}\sum_{\be,\mu}
g(\lnab{}dF(\be,\bmu),\lnab{}dF(\mu,\bbe))
=\frac{64\cos\theta}{\sin^{4}\theta}
\sum_{\be,\mu,\rho}(|g_{\be}\bmu\rho |^{2}+ 
|g_{\be}\bmu\brho |^{2})=}\non\\
&=& \frac{64\cos\theta}{\sin^{4}\theta}
\sum_{\be,\mu,\rho}( |g_{\be}\brho\mu |^{2} + 
|g_{\bmu}\be\brho |^{2})=\frac{128\cos\theta}
{\sin^{4}\theta}\sum_{\be,\mu,\rho}|g_{\be}\brho\mu|^{2}.
\end{eqnarray}
From (4.2) and (4.3),
\begin{eqnarray*}
&&\sum_{\be,\mu,\rho}Re\La{(}g_{\be}\mu\brho\, g_{\bbe}\rho\bmu\La{)}
= \sum_{\be,\mu}\sum_{\rho\neq \mu}|g_{\be}\brho\mu |^{2}
+\sum_{\be,\mu}Re\La{(}g_{\be}\mu\bmu\, g_{\bbe}\mu\bmu\La{)}\\
&&=\sum_{\be,\mu,\rho}|g_{\be}\brho\mu |^{2} -
\sum_{\be,\mu}|g_{\be}\bmu\mu |^{2}
+\sum_{\be,\mu}Re\La{(}g_{\be}\mu\bmu\, g_{\bbe}\mu\bmu\La{)}\\
&&=\sum_{\be,\mu,\rho}|g_{\be}\brho\mu |^{2}-\sum_{\be,\mu}
Re\La{(}\frac{i}{2}d\cos\theta(\be)g_{\bbe}\mu\bmu\La{)},
\end{eqnarray*}
so
\[(4.7)+(4.5)=
\frac{128\cos\theta}{\sin^{4}\theta}\sum_{\be,\mu}Re\La{(} \frac{i}{2}
d\cos\theta(\be)g_{\bbe}\mu\bmu \La{)}.\]
On the other hand, Proposition 3.2 and minimality of $F$ gives,
\[-\sum_{\be,\mu}\frac{4\cos\theta}{\sin^{2}\theta}
Re \La{(}ig_{\be}\bmu\mu\cdot \bbe\La{)}=
\frac{1-n}{4}\nabla\log\sin^{2}\theta = \frac{(n-1)\cos\theta}{2\sin^{2}\theta}\nabla \cos\theta.\]
Consequentely,
\begin{eqnarray*}
&&\frac{128\cos\theta}{\sin^{4}\theta}\sum_{\be,\mu}Re\La{(} \frac{i}{2}
d\cos\theta(\be)g_{\bbe}\mu\bmu \La{)}=
\frac{128\cos\theta}{\sin^{4}\theta}\sum_{\be,\mu}Re\La{(} -\frac{i}{2}
d\cos\theta(\bbe)g_{\be}\bmu\mu \La{)}\\[-1mm]
&&=-\frac{64\cos\theta}{\sin^{4}\theta}d\cos\theta \LA{(}Re\La{(} 
\sum_{\be,\mu}i g_{\be}\bmu\mu\cdot \bbe \La{)}\LA{)}=
\frac{8(n-1)\cos\theta}{\sin^{4}\theta}\|
\nabla \cos\theta\|^{2}.
\end{eqnarray*}
That is,
\begin{equation}
(4.7)+(4.5)=\frac{8(n-1)\cos\theta}{\sin^{4}\theta}\|
\nabla \cos\theta\|^{2}.
\end{equation}
 Using (3.9),
\begin{eqnarray}
\|\lnab{}\Jw\|^{2}&=&\sum_{\be}4 \langle \lnab{\be}\Jw\, ,\,\lnab{\bbe}
\Jw\rangle=\sum_{\be}\sum_{\mu,\rho}16\La{(} | \langle\lnab{\be}\Jw(\mu),
\rho\rangle |^{2} + |\langle\lnab{\be}\Jw(\bmu),\brho\rangle |^{2}\La{)}
\non\\[-2mm]
&=&64 \sum_{\be,\mu,\rho}\La{(}|\langle\lnab{\be}\mu,\rho\rangle|^{2} +  |\langle\lnab{\bbe}\mu,\rho\rangle|^{2}\La{)}.
\end{eqnarray}
Thus we see that $(4.6)=\frac{\cos\theta}{\sin^{2}\theta}\|\lnab{}
\Jw\|^{2}$. So we have  obtained the following formula:\\[-3mm]
\begin{Pp} If $N$ is K\"{a}hler-Einstein  with  Ricci tensor $Ricci^{N}=Rg$, and $F$ is a minimal immersion with equal K\"{a}hler angles, on an open set
without complex and Lagrangian points,\\[-10mm]
\begin{eqnarray}
\triangle\kappa &=& \cos\theta \La{(} -2nR
+\frac{32}{\sin^{2}\theta}\sum_{\be,\mu} R^{M}(\be,\mu,\bbe,\bmu)
 \non\\[-2mm]
&&~~~~~~~~~
+\frac{1}{\sin^{2}\theta}\|\lnab{}\Jw\|^{2}
+\frac{8(n-1)}{\sin^{4}\theta}\|\nabla\cos\theta\|^{2}~\La{)}.
\end{eqnarray}
\end{Pp}
Note that if $n=1$ we get the expression of Wolfson [W],~
$\triangle\kappa = -2R\cos\theta$.\\[-2mm]
\begin{Pp} If $N$ is K\"{a}hler-Einstein  with  Ricci tensor $Ricci^{N}=Rg$, and $F$ is a minimal imersion with equal K\"{a}hler angles, then: \\[1mm]
(i)~ If $n=2$,\\[-9mm]
\begin{equation}
\int_{M}\!nR\sin^2\theta\cos^2\theta \, Vol_{M} = 0.
\end{equation} 
(ii)~~  If $n\geq 3$ and $F$ has no complex points,\\[-4mm]
\begin{equation}
\int_{M}\!nR\sin^2\theta\cos^2\theta \, Vol_{M} = \int_{M}(n-2)(n-2+2 \cot^{2}\theta)\|\nabla\cos\theta\|^{2}\,Vol_{M}.\\[3mm]
\end{equation}
\end{Pp}
\em Proof. \em 
Multiplying (4.10) by $\sin^2\theta\cos\theta$, we get, on $M\sim {\cal C}
\cup {\cal L}$, and using Lemma 3.2,
\begin{eqnarray*}
\sin^2\theta\cos\theta\triangle\kappa &=& -2n\sin^2\theta\cos^2\theta R
+ 2\langle S\Fw,\Fw\rangle\\[-1mm]
&&+\cos^{2}\theta\|\lnab{}\Jw\|^{2}
+\frac{8(n-1)\cos^{2}\theta}{\sin^{2}\theta}\|\nabla\cos\theta\|^{2}.
\end{eqnarray*}
On the other hand, $\kappa = n\log \left( \frac{1+\cos\theta}{1-\cos\theta}\right)$, and so, 
$\triangle\kappa =\frac{2n}{\sin^{2}\theta}
\triangle\cos\theta +\frac{4n\cos\theta}{\sin^{4}\theta}
\|\nabla\cos\theta\|^{2}$. Hence,\\[-8mm]
\begin{eqnarray}
\lefteqn{2n\cos\theta\triangle\cos\theta +\frac{4n\cos^{2}\theta}{\sin^{2}\theta}
\|\nabla\cos\theta\|^{2}=}\\[-2mm]
&\!\!\!\!\!\!=&\!\!\!-2n\sin^2\theta\cos^2\theta R
+ 2\langle S\Fw,\Fw \rangle +\cos^{2}\theta\|\lnab{}\Jw\|^{2}
+\frac{8(n\!-\!1)\cos^{2}\theta}{\sin^{2}\theta}\|\nabla\cos\theta\|^{2}.
\non
\end{eqnarray}
Recall that, from (3.1), and considering $\Fw$ a 2-form,
 $\|\lnab{}\Fw\|^{2}=\frac{1}{2}\cos^2\theta\|\lnab{}\Jw\|^{2}+
n\|\nabla\cos\theta\|^{2}$.
Since $\triangle\cos^{2}\theta = 2\cos\theta\triangle\cos\theta +
2\|\nabla\cos\theta\|^{2}$, substituting this into (4.13),
 we have\\[-6mm]
\begin{equation}
n\triangle\cos^{2}\theta = -2n\sin^2\theta\cos^2\theta R
+ 2\langle S\Fw,\Fw \rangle 
+2\|\lnab{}\Fw\|^{2}
+\frac{4(n\!-\!2)\cos^{2}\theta}{\sin^{2}\theta}\|\nabla\cos\theta\|^{2}
\end{equation}\\[-8mm]
and, for $n=2$,\\[-8mm]
\begin{eqnarray}
n\triangle\cos^{2}\theta &=& -2n\sin^2\theta\cos^2\theta R
+ 2\langle S\Fw,\Fw \rangle +2\|\lnab{}\Fw\|^{2}.
\end{eqnarray}
Let us now suppose that $n\geq 3$. 
Then, under the condition of no complex points, (4.14) is valied on $\Omega^{0}_{2n}$ and  also on $\Omega^{0}_{0}$. 
From smoothness over all $M$ of all maps into consideration 
(the first three terms 
of the right-hand side of (4.14) are smooth, and the last term is also 
smooth for $n\neq 2$), and the fact that 
the set $M\sim (\Omega^{0}_{0}\cup\Omega^{2n}_{0})$ is a set of 
Lagrangian points with no interior, 
 formula (4.14) is valid on all $M$.
Integrating over $M$, and using (3.17), we have
\[
\int_{M}\!2nR\sin^2\theta\cos^2\theta \, Vol_{M}\! =\!\!
\int_{M}\!\!\LA{(}\!\! -2(n\!-\!(n-2)^2) +\frac{4(n\!-\!2)\cos^{2}\theta}{\sin^{2}\theta}
+2n\LA{)}\|\nabla\cos\theta\|^{2}\,Vol_{M}, \]
leading to (4.12). If $n=2$, we see that (4.15)
 is also valid at Lagrangian and complex points. In fact (see 
Lemma 3.2 and (3.1)), 
all terms of (4.15) vanish at interior points of the Lagrangian
and complex sets. Since they are smooth on all $M$, they must vanish 
at boundary points of its complementary in $M$. Thus, the above equation
is valid on all $M$, with or without complex or
Lagrangian points, and  all its terms are smooth. Then, (4.11) follows
by integration on $M$ of (4.15), and use of (3.17).
~~~~~~~\qed\\[6mm]
\em Proof of Theorem 1.2. and  Theorem 1.3 \em ~~ If $n=2$
and $R\neq 0$, $(4.11)$
implies $\sin^2\theta\cos^2\theta =0$. Hence  $F$ is either Lagrangian
or a complex submanifold.
If $n\geq 3$, and $F$ has no complex points, 
 the right-hand side of (4.12) is non-negative, while the
left-hand side is non-positive for $R<0$. Then,
 $\sin^2\theta\cos^2\theta =0$
must hold on all $M$, that is, $F$ is Lagrangian.
If  $R=0$, the right-hand side of (4.12) must vanish.
Then, for $n\geq 3$, 
$\cos\theta$ must be constant, and we have proved Theorem 1.2. 
If $\cos\theta$ is constant, and if $F$ is not a complex submanifold,
the right-hand side of (4.12) vanishes. Hence, 
if $R\neq 0$, $F$ is Lagrangian, and Theorem 1.3 is proved.~~~~~\qed
\\[5mm]
\em Proof of Theorem 1.4. \em   If $M$ is not Lagrangian, under the curvature
condition on $M$, by
Proposition 3.4, for $n=2$, or $3$,  $(M,\Jw,g_{M})$
is a K\"{a}hler manifold and 
 $\cos\theta$ is constant. So, if $M$ is not a complex submanifold,
it has no complex directions, and
by (4.11), or (4.12),  $R=0$. In general, if $n\geq 1$ and $\theta$ is constant, Proposition 3.4 also applies.~~\qed
\\[5mm]
Under the conditions of Theorem 1.4, if $M$ is homeomorphic to a 4 or a 6 dimensional sphere, minimaly immersed into a K\"{a}hler-Einstein manifold,
and with equal K\"{a}hler angles, 
then  it must be  Lagrangian, for it is well known that such manifolds 
cannot carry a K\"{a}hler structure. Obviously, any Riemannian 
manifold $M$ with strictly positive isotropic scalar curvature cannot carry any K\"{a}hler
structure. Moreover, such condition for $n=2$ would imply $M$ to
be homeomorphic to a 4-sphere. We also remark that we only need to
require $S_{isot}(\{Z_{\al}\})\geq 0$ on the maximal totally isotropic
subspace $\{Z_{\al}\}$ defined by a  diagonalizing orthonormal basis
of $\Fw$, and outside Lagrangian points, to obtain the same conclusion
given in Theorem 1.4.

As an observation, Theorem 1.4 should be compared with the following lemma:
\begin{Lm} Let $F$ be a minimal immersion, and $n\geq 2$.
 If $\cos\theta$ is constant $\neq 1,0$, then\\[2mm]
$(i)$~~$(A,B,C)\ra \gf{A}{B}{C}$ is symmetric whenever $A,B$, and $C$ are
not all of the same type.\\[1mm]
$(ii)$~~$\langle\lnab{\bbe}\mu,\ga\rangle =0,~~~\forall \be,\mu,\ga$.\\[1mm]
$(iii)$~ $\Fw$ is an harmonic 2-form.\\[1mm]
$(iv)$~$32\sum_{\be,\mu}R^{M}(\be,\mu,\bbe,\bmu) = -64\sum_{\be,\mu,\rho}|\langle \lnab{\be}\mu,\rho\rangle |^{2}=-\|\lnab{}\Jw\|^{2}\leq 0.$
\end{Lm}
\em Proof. \em  Since $\cos\theta$ is constant, we obtain
$(4.3)=0$. This, together (4.2),  and the symmetry of $\lnab{}dF$,
proves $(i)$. But $(i)$ and (4.1) implys $(ii)$. $(iii)$ comes from
(3.16).
Now we prove $(vi)$. Since $\Fw$ is harmonic, from Weitzenb\"{o}ck formula
(3.14) we conclude $\langle S\Fw, \Fw\rangle=-\|\lnab{}\Fw\|^{2}$.
 Lemma 3.2 and (3.1) (but considering $\Fw$ a 2-form)
gives  $(iv)$.~~~~~~\qed
\\[5mm]
\em Remark 3. \em  If $N$ is a K\"{a}hler manifold of constant 
 holomorphic sectional curvature equal to $K$ ( and so
$R=\frac{(2n+1)K}{2}$), and the isotropic scalar curvature of $M$ satisfies $S_{isot}\geq c= constant$, we get  
from Gauss equation, with  $\{X_{\al},Y_{\al}\}$  a diagonalizing
orthonormal basis of $\Fw$,
\begin{equation}
\sum_{\rho,\mu}R^{M}(\mu,\rho,\bmu,\brho)=\frac{n(n-1)}{16}\sin^{2}\theta K
-\sum_{\rho,\mu}\|\lnab{}dF(\mu,\brho)\|^{2},
\end{equation}
that $c\leq \frac{n(n-1)K}{4}$.
Thus,  non-negative isotropic scalar curvature on $M$ is a possible condition for $K \geq 0$. In the case $K=0$, that is, $N$ is the flat complex torus, 
then (4.16) (with $K=0$) is valied for any orthonormal basis
$\{X_{\al},Y_{\al}\}$. This implies that, 
 for $n\geq 2$,  $F$ must be  totally geodesic, and so $M$ is flat.

We also note that if $c=\frac{nR}{4}$,
 the right-hand side of (4.10) becomes $>0$, outside  Lagrangian points. An
application  of the maximum principle at a maximum point of $\kappa$ would
conclude that $F$ must be Lagrangian. But such a
 lower bound  $c$ is not  possible for
the scalar isotropic curvature of $M$ minimaly immersed in $N$ with
constant   holomorphic sectional curvature  $K>0$. 
\\[6mm]
\em Remark 3. \em  If $n\geq 2$ and $F$ is a pluriminimal immersion with equal K\"{a}hler angles  into a K\"{a}hler-Einstein  manifold $N$,
and $F$ is not a complex submanifold, then $N$ must be  Ricci-flat. Moreover, since $F$  has cons\-tant equal K\"{a}hler angles,
 the scalar isotropic curvature of $M$ with respect to
 the maxi\-mal isotropic subspace defined by a diagonalizing orthonormal basis
of $\Fw$
 will be $\leq 0$, with equality to zero iff 
$(M,\Jw,g_{M})$ is K\"{a}hler (see Lemma 4.1). We leave the following
question: Is $(M,\Jw,g_{M})$ K\"{a}hler manifold a sufficient condition for a minimal immersion $F$, with constant equal K\"{a}hler angle,  immersed  into a Ricci-flat K\"{a}hler manifold $N$,
 to be pluriminimal?  
If $N$ is the flat complex torus
and $F:M\ra N$ is minimal, under the conditions stated in the question, 
the Gauss equation   implies that  $F$ is pluriminimal. 
A way to find  plu\-ri\-mi\-ni\-mal submanifolds in hyper-K\"{a}hler
manifolds is given in the next example, where the assumption of
non-negative isotropic curvature does not imply necessarely
$F$ totally geodesic (and $M$ flat), since hyper-K\"{a}hler
manifolds do not need to be flat.\\[6mm]
\em Example. \em  Let $(N,I,J,g)$ be an hyper-K\"{a}hler manifold of real dimension $8$. Thus,
$I$ and $J$ are two $g$-orthogonal complex structures on $N$, such that 
$IJ=-JI$ and $\lnab{}I=\lnab{}J=0$, where $\lnab{}$ is the Levi-Civita connection relative to $g$. It is known that such manifolds are Ricci-flat 
([B]). Set $K=IJ$. For each $\nu$, $\phi$, 
we take $``\nu\phi"=(\cos\nu, \sin\nu\cos\phi,\sin\nu\sin\phi)\in S^{2}$, and define  $J_{\nu\phi}=\cos\nu I+\sin\nu\cos\phi J
+\sin\nu\sin\phi K$. These $J_{\nu\phi}$ are the   complex structures
on $N$ compatible with its hyper-K\"{a}hler structure, that is, 
 they are $g$-orthogonal and  $\lnab{}J_{\nu\phi}=0$. 

Two of such complex structures, $J_{\nu\phi}$ and $J_{\mu\rho}$, anti-commute
at a point $p$  iff  $ J_{\nu\phi}(X)$ and $J_{\mu\rho}(X)$
are orthogonal for some non-zero $X\in T_{p}N$, iff 
$\nu\phi$ and $\mu\rho$ are orthogonal in $\R{3}$.
Thus,  they anti-commute at a point $p$ iff they anti-commute everywhere.
If that is the case $J_{\nu\phi}\circ J_{\mu\rho}= J_{\sigma\epsilon}$,
where  $\{\nu\phi, \mu\rho,\sigma\epsilon\}$ is a direct orthonormal
basis of $\R{3}$.
For each  unit vector $X\in T_{p}N$, set 
$H_{X}=span\{ X,IX,JX,KX\}$ $=span\{X,J_{\nu\phi}(X),
J_{\mu\rho}(X), J_{\sigma\epsilon}(X)\}$, for any orthonormal basis $\{
\nu\phi, \mu\rho,\sigma\epsilon\}$.
If $Y\in H_{X}^{\bot}$ is another unit vector,
 then  $H_{X}\bot H_{Y}$.
Let $\omega_{\nu\phi}$ be the K\"{a}hler form
of $(N,J_{\nu\phi},g)$. 
Let $E$ be a 4-dimensional vector sub-space of $T_{p}N$. 
 We first note that  $E=H_X$ for some $X\in E$, iff $J_{\nu\phi}(E)\subset E$ for any $\nu,\phi$. If that is the case, then $E$ is not a Lagrangian subspace with respect to any complex structure $J_{\mu\rho}$.
In general,  $E$ contains a $J_{\nu\phi}$-complex line for some $\nu\phi$ 
iff $dim(E\cap H_{X})\geq 2$ for some $X\in E$.  If that is the case, and if
$E$ is a Lagrangian subspace of $T_{p}N$ with respect to $J_{\mu\rho}$,
then $ \nu\phi\bot\mu\rho$. Furthermore, if $E$ is a 
$J_{\nu\phi}$-complex subspace, then $E$ is $J_{\mu\rho}$-Lagrangian iff
there exist an orthonormal basis $\{X,J_{\nu\phi}X,Y,J_{\nu\phi}Y\}$ of $E$ 
with $H_X\bot H_Y$. To see this, let us suppose $E$ is
$J_{\nu\phi}$-complex subspace and
$J_{\mu\rho}$-Lagrangian.  We take $\{X, J_{\nu\phi}X, Y, J_{\nu\phi}Y\}$
an ortonormal basis of $E$. Then 
$Y\in span\{X,J_{\nu\phi}X,J_{\mu\rho}X\}^{\bot}$. So $Y=t
J_{\sigma\epsilon}X + \tilde{Y}$, for some
$t\in\R{}$ and $\tilde{Y}\in H_{X}^{\bot}$, and 
where $\{\nu\phi, \mu\rho, \sigma\epsilon\}$ is an ortonormal basis of 
$\R{3}$. As $E\neq H_{X}$, $\tilde{Y}\neq 0$. From 
$0=\langle J_{\mu\rho}Y, J_{\nu\phi}X\rangle$, we get $t=0$. Thus, $Y\in H_{X}^{\bot}$. We observe that, in general, 
 $J_{\mu\rho}$-Lagrangian subspaces
do not need to be $J_{\nu\phi}$-complex, as for example $E=\{ X,J_{\nu\phi}X,
Y,J_{\sigma\epsilon}Y\}$, with $Y\in H_{X}^{\bot}$, that contains two 
orthogonal complex lines for different complex strutures.

Any $J_{\nu\phi}$-complex submanifold
$F:M\ra N$ of real dimension  $4$, such that, for each point $p\in M$, there exist an orthonormal basis $\{X,J_{\nu\phi}X,Y,J_{\nu\phi}Y\}$ of $T_{p}M$ 
with $H_X\bot H_Y$, is, for each $\mu\rho$,
 a minimal submanifold of $(N,J_{\mu\rho},g)$  with 
constant equal K\"{a}hler angles,
 and  $\pm J_{\nu\phi}$ is also the complex
structure of $M$ which comes from polar decomposition of $\omega_{\mu\rho}$ 
restricted to $M$. 
In fact, such an orthonormal basis of $T_{p}M$ diagonalizes $\omega_{\mu\rho}$ restricted to $M$, and the K\"{a}hler angle $\theta$ is such that 
$\cos\theta=\pm \langle \nu\phi, \mu\rho\rangle$, where $<,>$ is the inner product of $\R{3}$.
Next proposition is an application of Theorem 1.4, for 4-dimensional submanifolds of $N$, where   $\omega_{I}$ is the K\"{a}hler form of $(N,I,g)$:
\begin{Pp} Let $F:M\ra N$ be a minimal immersion of a compact, oriented  
4-dimensional sub\-ma\-ni\-fold  with non-negative isotropic scalar 
curvature, and  such that $\forall \nu\phi\in S^{2}$, $F$
 has equal K\"{a}hler angles with respect to  $J_{\nu\phi}$.
If $\exists p\in M$ and  $\exists X\in T_{p}M$, unit vector,
such that $dim(T_{p}M\cap H_X)\geq 2$, then there exists 
$\nu\phi\in S^{2}$ such that $M$ is a $J_{\nu\phi}$-complex 
submanifold. Furthermore, if $J_{\nu\phi}=I$ then $F: M \ra (N,I,g)$ is obviously pluriminimal. If $J_{\nu\phi}\neq I$ but
 $T_{p}M\cap H_X^{\bot}\neq \{ 0\}$, 
then $F^{*}\omega_I =\cos\nu J_{\nu\phi}$, and if  $F$ is not 
$J_{I}$-Lagrangian, 
 $F: M \ra (N,I,g)$ is still pluriminimal.
\end{Pp} 
Note that, if $T_{p}M=H_{X}$, then $J_{\nu\phi}$ can be chosen equal to $I$.
The first conclusion of this result is the 4-dimensional version of a result of 
Wolfson [W], for  $M$ a real surface and $N$ a Ricci-flat K3 surface. 
In the latter case, there is only one K\"{a}hler angle,
$\forall X$ $dim(T_{p}M\cap H_X)=2$  is automatically satisfied,
and the isotropic scalar curvature is always zero.\\[2mm]
\em Proof. \em  From the  assumption, $dim(T_{p}M\cap H_X)\geq 2$, 
 we may take a unit vector $Z\in T_{p}M\cap H_X$ such that  $Z\bot X$.
 Then, $Z=J_{\nu\phi}(X)$ for some $\nu\phi$.
Thus, $span\{X,J_{\nu\phi}(X)\}\subset T_p M$. This implies 
$F^{*}\omega_{\nu\phi}(X,J_{\nu\phi}(X))=1$. As the
K\"{a}hler angles are equal, $\cos\theta_{\nu\phi}=1$ at $p$.
Applying Theorem 1.4  to $F:M \ra (N,J_{\nu\phi},g)$, 
$  F^{*}\omega_{\nu\phi}= \cos \theta_{\nu\phi}J_{\omega_{\nu\phi}}$
with $\cos \theta_{\nu\phi}$ constant. 
Then $\cos\theta_{\nu\phi}=1$ everywhere. That is, $M$ is a $J_{\nu\phi}$-complex submanifold. 
Moreover, from the second assumption, $T_{p}M\cap H_X^{\bot}\neq \{0\}$,  
we may take a unit vector $Y\in T_{p}M\cap H^{\bot}_{X}$. Then $\{X,J_{\nu\phi}X,
Y,J_{\nu\phi}Y\}$ constitutes an orthonormal  basis of $T_{p}M$,
that diagonalizes $F^{*}\omega_{I}$, and  $F^{*}\omega_{I}=\cos\nu J_{\nu\phi}$. This means that
$\nu$ or $\nu+\pi$ is the constant K\"{a}hler angle of $F:M\ra (N,I,g)$,
and, since $M$ is a $J_{\nu\phi}$-complex submanifold, it is pluriharmonic
with respect to $\pm J_{\nu\phi}$, and so, if $\cos\nu\neq 0$,
it is pluriminimal as an immersion into $(N,I,g)$.~~~~~~\qed
\section{Appendix: The computation of $\triangle \kappa$}
\setcounter{Th}{0}
\setcounter{Pp}{0}
\setcounter{Cr} {0}
\setcounter{Lm} {0}
\setcounter{Def} {0}
\setcounter{equation} {0}
We prove (1.6) for $F$ minimal and outside complex and Lagrangian points. 
First, we compute some derivative formulas of a determinant, 
which we will need.
\begin{Lm} Let $A:M\ra {\cal M}_{m\times m}(\Co)$ be a smooth map of
matrices $ p\ra A(p)=[A_{1},\ldots, A_{m}]$, where $A_{i}(p)$ is a column
vector of $~\Co^{m}$ and $M$ is a Riemannian manifold with its Levi-Civita
connection $\lnab{}$. Assume that, at $p_{0}$, $A(p_{0})$ is a
diagonal matrix $D=D(\lambda_{1},\ldots, \lambda_{m})$. Then, at $p_{0}$
\\[-5mm]
\[ d~(\mbox{det~}A)(Z) = \sum_{1\leq j\leq m} \La{(}\prod_{k\neq
j}\lambda_{k}\La{)} dA_{j}^{j}(Z),\]
\\[-12mm]
\begin{eqnarray*}
\lefteqn{\mbox{Hess~}(det~A)(Z,W)=\lnab{} d(det A)(Z,W)=}\\ &=&\sum_{1\leq
j,k \leq m} \La{(}\prod_{s\neq j,k}\lambda_{s}\La{)}\mbox{det} \left[
\begin{array}{cc}
dA_{j}^{j}(Z) & dA^{k}_{j}(Z)\\
dA_{k}^{j}(W) & dA_{k}^{k}(W)
\end{array}\right] +\sum_{1\leq j\leq m} \La{(}\prod_{s\neq
j}\lambda_{s}\La{)} \mbox{Hess~}A_{j}^{j}(Z,W).\end{eqnarray*} 
In particular,
if $e_{1},\ldots,e_{r}$ is an orthonormal  basis of $T_{p_{0}}M$, then, at
$p_{0}$,\\[-9mm]
\begin{eqnarray*}
\lefteqn{\triangle (det~A)= Trace~\mbox{Hess~}(det ~A)=}\\[-2mm]
&=&\sum_{1\leq\al\leq r}~~\sum_{1\leq j,k \leq m} \La{(}\prod_{s\neq
j,k}\lambda_{s}\La{)}\mbox{det}
\left[ \begin{array}{cc}
dA_{j}^{j}(e_{\al}) & dA^{k}_{j}(e_{\al})\\ dA_{k}^{j}(e_{\al}) &
dA_{k}^{k}(e_{\al}) \end{array}\right] +\sum_{1\leq j\leq m}
\La{(}\prod_{s\neq j}\lambda_{s}\La{)} \triangle A_{j}^{j}.
\end{eqnarray*}
\end{Lm}
On each $\Omega_{2k}^{0}$, the 
complex structure $\Jw$ and  the  sub-vector bundle
${\cal K}_{\omega}^{\bot}$ are smooth.
Moreover, $J_{\omega}$ is $g_{M}$-orthogonal.
Thus,  for each $p_{0}\in \Omega_{2k}^{0}$, there exists a locally
$g_{M}$-orthonormal  frame of ${\cal K}_{\omega}^{\bot}$ defined 
on a neighbourhood of  $p_{0}$, of the
form $X_{1},J_{\omega}X_{1},\ldots,X_{k},J_{\omega}X_{k}$. We enlarge
this frame to a $g_{M}$-orthonormal local frame on $M$, on a 
neighbourhood of $p_{0}$:
\begin{equation}
X_{1}\;, Y_{1}=J_{\omega}X_{1}\;,\ldots,X_{k}\;,Y_{k}=J_{\omega}X_{k}\;,
\,X_{k+1}\;,Y_{k+1}\;,\ldots\;,
X_{n}\;,Y_{n}\end{equation}
where $X_{k+1},Y_{k+1},\ldots X_{n},Y_{n}$ is any $g_{M}$-orthonormal  
frame of
${\cal K}_{\omega}$, and which at $p_{0}$ is a
dia\-go\-na\-li\-zing basis of $F^{*}\omega$. Note that in general it is not
possible to get smooth diagonalizing $g_{M}$-orthonormal  frames in a whole
neighbourhood of a point $p_{0}$, unless , for instance, $\Fw$ has equal K\"{a}hler angles. We use the notations in section 3.1.
We define a local complex structure
on a neighbourhood of $p_{0}\in \Omega^{0}_{2k}$ as
$\tilde{J}=J_{\omega}\oplus J'$, where $J_{\omega}$ is defined only on
${\cal K}^{\bot}_{\omega}$, and $J'$ is the local complex structure on
${\cal K}_{\omega}$, defined on a neighbourhood of $p_{0}$ by\\[-5mm]
\begin{equation}
J'Z_{\al}=iZ_{\al},~~~~J'Z_{\bal}=-iZ_{\bal}, ~~\forall \al\geq k+1.
\end{equation}\\[-5mm]
Thus, the vectors  $Z_{\al}$ are of type (1,0) with respect to
$\tilde{J}$, for $\forall \al$.
Since $\tilde{J}$ is $g_{M}$-orthogonal, then, $\forall \al,\be$, on a neighbourhood of $p_{0}$,\\[-4mm] 
\begin{equation}
\langle \lnab{Z}\tilde{J}({\al}),{\be}\rangle =2i\langle\lnab{Z} {\al},
{\be}\rangle= -\langle
{\al},\lnab{Z}\tilde{J}({\be})\rangle,~~~~~~~~~~ \langle
\lnab{Z}\tilde{J}({\al}),{\bbe}\rangle =0, 
\end{equation}\\[-5mm]
Note that $\Fw$ and $\tg$, where  $\tg$  is given in (1.1),
as 2-tensors, are both of type $(1,1)$  with respect to $\tilde{J}$, and have the same kernel $\kw$. They are related by $\tg(X,Y)=\Fw(X,\Jw Y)=\Fw(X,\tilde{J}Y)$.
Set $\tg_{AB}=\tg({A},{B})$, 
and define $\overline{\bar{B}}=B$,  $\forall A,B\in\{1,\ldots,n, \bar{1},\ldots, \bar{n}\}$, and  set $\epsilon_{\al}=+1$, $\epsilon_{\bal}=-1$, 
 $\forall 1\leq \al\leq n$. Let $ 1\leq \al, \be \leq n$,
$ A,B\in\{1,\ldots,n, \bar{1},\ldots, \bar{n}\}$,
 and $C\in\{1,\ldots,n\}\cup \{\overline{k+1},\ldots,\bar{n}\}$. Then
\begin{equation}\left.\begin{array}{ll}
F^{*}\omega({\al},{C})=g(JdF({\al}),dF({C}))=0 &~~~~~~\forall p
\mbox{~near~} p_{0}\\[1mm]
F^{*}\omega({\al},{\bbe})=g(JdF({\al}),dF({\bbe}))
=\frac{i}{2}\delta_{\al\be}\cos\theta_{\al} &~~~~~~\mbox{at~}p_{0}\\[1mm]
\tg_{AB}=i\epsilon_{B}F^{*}\omega({A},{B})=i\epsilon_{B}g(JdF({A}),dF({B}))
&~~~~~~\forall p \mbox{~near~} p_{0}\\[1mm] 
\tg_{\al C}=\tg_{\bal\bar{C}}=0
&~~~~~~\forall p \mbox{~near~} p_{0}\\[1mm]
\tg_{\al\bbe}=\tg_{\bal\be}=\frac{1}{2}\delta_{\al\be}\cos\theta_{\al}
&~~~~~~\mb
{at~}p_{0} \end{array}\right\} . \end{equation}
At a point $p_{0}$, with K\"{a}hler angles $\theta_{\al}$,
$g_{M}\pm\tg$ is represented  in the
unitary  basis $\{\sqrt{2}\al,\sqrt{2}\bal\}$, 
by the diagonal matrix $ g_{M}\pm\tg=D(1\pm\cos\theta_{1},\ldots,1\pm\cos\theta_{n},
1\pm\cos\theta_{1},\ldots,1\pm\cos\theta_{n})$,
and so
\begin{equation}
 det (g_{M}\pm\tg)=\prod_{1\leq\al\leq n}(1\pm\cos\theta_{\al})^{2}.
\end{equation}
If $p_{0}$ is a point without complex directions, $\cos\theta_{\al}\neq 1$,
$\forall\al\in\{1,\ldots,n\}$, then $\tg<g_{M}$. 
Thus, on a neighbourwood of $p_{0}$, we may consider the map $\kappa$.
 \begin{equation}
\kappa =\frac{1}{2}\log \left( \frac{det(g_{M}+\tg)}{det(g_{M}-\tg)}\right)
=\sum_{1\leq \al\leq
n}\log\LA{(}\frac{1+\cos\theta_{\al}}{1-\cos\theta_{\al}}\LA{)}. 
\end{equation}
This map is continuous outside the complex points, and smooth on
each $\Omega_{2k}^{0}$.
 We wish to compute $\triangle \kappa$ on $\Omega_{2k}^{0}$. 
\begin{Lm} At $p_{0}\in \Omega_{2k}^{0}$, without complex directions and 
for $Z, W\in T_{p_{0}}M$,
\[d(det(g_{M}\pm \tg))(Z)= \pm 4\sum_{1\leq \mu\leq n}\frac{
\prod_{1\leq\al\leq n}(1\pm\cos\theta_{\al})^{2}}{(1\pm\cos\theta_{\mu})}
d\tg_{\mu\bmu}(Z),\]\\[-12mm]
\begin{eqnarray*}
\lefteqn{Hess\la{(}det(g_{M}\pm \tg)\la{)}(Z,W)=}\\
&=&16 \La{(}\prod_{1\leq\al\leq n}(1\pm\cos\theta_{\al})^{2} \La{)}
\sum_{\mu,\rho}\frac{1}{(1\pm\cos\theta_{\mu})(1\pm\cos\theta_{\rho})}
d\tg_{\mu\bmu}(Z)d\tg_{\rho\brho}(W)\\[-2mm]
&& -8\La{(}\prod_{1\leq\al\leq n}(1\pm\cos\theta_{\al})^{2} \La{)}
\sum_{\mu,\rho}\frac{1}{(1\pm\cos\theta_{\mu})(1\pm\cos\theta_{\rho})}
d\tg_{\mu\brho}(W)d\tg_{\rho\bmu}(Z)\\[-2mm]
&&\pm 4 \La{(}\prod_{1\leq\al\leq n}(1\pm\cos\theta_{\al})^{2} \La{)}
\sum_{\mu}\frac{1}{(1\pm\cos\theta_{\mu})} Hess\tg_{\mu\bmu}(Z,W).
\end{eqnarray*}
\end{Lm}
\em Proof. \em
Using the unitary basis $\{\sqrt{2}{\al},\sqrt{2}{\bal}\}$ of $T_{p}^{c}M$,
for $p$ near $p_{0}$, $g_{M}+\tg$ is represented by the matrix 
\[
g_{M}\pm\tg = \left[\begin{array}{cc}
g_{M}\pm\tg(\sqrt{2}\al,\sqrt{2}\bga) & g_{M}\pm\tg(\sqrt{2}\al,\sqrt{2}\ga)\\
g_{M}\pm\tg(\sqrt{2}\bal,\sqrt{2}\bga) & g_{M}\pm\tg(\sqrt{2}\bal,\sqrt{2}\ga)
\end{array}\right]=
\left[\begin{array}{cc}
\delta_{\al\ga}\pm 2\tg_{\al\bga}& 0 \\
0 & \delta_{\al\ga}\pm 2\tg_{\bal\ga}\end{array}\right] 
\]
that at $p_{0}$ is the diagonal matrix $
D(1\pm \cos\theta_{1},\ldots,1\pm \cos\theta_{n}, 1 \pm
\cos\theta_{1},\ldots,1\pm \cos\theta_{n})$. The lemma follows as a simple
application of lemma 5.1, and noting that
$\tg_{\mu\brho}=\tg_{\brho\mu}$.~~~~~\qed\\[4mm] 
On $\Omega_{2k}^{0}$,\\[-8mm] 
\begin{eqnarray*}
2\triangle \kappa\!\!\!
&=&\!\! \triangle \log (det (g_{M}+\tg)) -\triangle \log (det(g_{M}-\tg))\\
\!\!\!\!&=&\!\!\!\!\frac{\triangle (det (g_{M}+\tg))}{det(g_{M}+\tg)}
-\frac{\|d(det (g_{M}+\tg))\|^{2}}{(det (g_{M}+\tg))^{2}} -\frac{\triangle
(det (g_{M}-\tg))}{det(g_{M}-\tg)} +\frac{\|d(det (g_{M}-\tg))\|^{2}}{(det
(g_{M}-\tg))^{2}} .\\[-8mm]
\end{eqnarray*}
From the above lemma and\\[-8mm]
\begin{eqnarray*}
\|d(det (g_{M}\pm\tg))\|^{2}&=&
4\sum_{\be}d(det (g_{M}\pm\tg))({\be})d(det(g_{M}\pm\tg))({\bbe}) \\[-2mm]
\triangle det(g_{M}\pm \tg)
&=& 4\sum_{\be}Hess(det(g_{M}\pm\tg))({\be},{\bbe}) 
\end{eqnarray*}\\[-10mm]
we have at $p_{0}$,
\begin{equation}
2\triangle \kappa = \sum_{\be,\mu,\rho}\frac{64(\cos\theta_{\mu}+
\cos\theta_{\rho})}{\sin^{2}\theta_{\mu}\sin^{2}\theta_{\rho}}
d\tg_{\mu\brho}({\bbe})d\tg_{\rho\bmu}({\be}) +
\sum_{\be,\mu}\frac{32}{\sin^{2}\theta_{\mu}}
Hess\tg_{\mu\bmu}({\be},{\bbe}).
\end{equation}
Recalling (2.4), and  $d(\Fw(X,Y))(Z)=\lnab{Z}\Fw(X,Y)+\Fw(\lnab{Z}X,Y)
+\Fw(X,\lnab{Z}Y)$, using (5.4), we obtain
\begin{Lm}
$\forall p$ near $p_{0}\in\Omega^{0}_{2k}$, $Z\in T^{c}_{p}M$, and
$\mu,\ga\in\{1,\ldots,n\}$\\[-8mm]
\begin{eqnarray*}
d\tg_{\mu\bga}(Z) & = & i\gf{Z}{\mu}{\bga}-i\gf{Z}{\bga}{\mu} 
 +2\sum_{ \rho}\La{(}\langle
\lnab{Z}{\mu},{\brho}\rangle\tg_{\rho\bga} +\langle
\lnab{Z}{\bga},{\rho}\rangle\tg_{\mu\brho}\La{)}\\[-3mm] 
0=d\tg_{\mu\ga}(Z)& =& -i\gf{Z}{\mu}{\gamma}+i\gf{Z}{\ga}{\mu}
 +2\sum_{\rho}\La{(}\langle
\lnab{Z}{\mu},{\rho}\rangle\tg_{\brho\ga} -\langle
\lnab{Z}{\ga},{\rho}\rangle\tg_{\mu\brho}\La{)}.
 \end{eqnarray*}\\[-8mm]
In particular, at $p_{0}$\\[-10mm]
\begin{eqnarray*}
d\tg_{\mu\bga}(Z) & = & i\gf{Z}{\mu}{\bar{\gamma}}
-i\gf{Z}{\bga}{\mu}-(\cos\theta_{\mu} -\cos\theta_{\ga})\langle
\lnab{Z}{\mu},{\bga}\rangle\\[-1mm]
0=d\tg_{\mu\ga}(Z)	& = & -i\gf{Z}{\mu}{\gamma}
+i\gf{Z}{\ga}{\mu}+(\cos\theta_{\mu} +\cos\theta_{\ga})\langle \lnab{Z}{\mu},{\ga}\rangle.
\end{eqnarray*}
\end{Lm}
\begin{Lm} If F is minimal and $p_{0}\in \Omega^{0}_{2k}$ is a point without 
complex directions, then for each $\mu\in\{1,\ldots,n\}$
\begin{eqnarray*}
\sum_{1\leq\be\leq n}\lefteqn{Hess\tg_{\mu\bmu}(\be,\bbe)=
\sum_{1\leq\be\leq n}d\La{(}d\tg_{\mu\bmu}
(\be)\La{)}(\bbe)-d\tg_{\mu\bmu}(\lnab{\bbe}\be)
=} \\
~~= \sum_{1\leq\be\leq n}&& i R^{N}(dF(\be),dF(\bbe),dF(\mu),
JdF(\bmu)+i\cos\theta_{\mu}dF(\bmu)) \\[-3mm]
&&+2 Im \La{(}R^{N}(dF(\be),dF(\mu),dF(\bbe),
JdF(\bmu)+i\cos\theta_{\mu}dF(\bmu))\La{)} \\
&&+2\sum_{1\leq \rho\leq n}\frac{(\cos\theta_{\rho}-\cos\theta_{\mu})}
{\sin^{2}\theta_{\rho}}\;
\La{(}|\gf{\be}{\mu}{\rho}|^{2}+|\gf{\be}{\bmu}{\brho}|^{2}\La{)}\\[-2mm]
&&-2\sum_{1\leq \rho\leq n}\frac{(\cos\theta_{\rho}+
\cos\theta_{\mu})}{\sin^{2}\theta_{\rho}}\;\La{(}
|\gf{\be}{\mu}{\brho}|^{2}+|\gf{\be}{\bmu}{\rho}|^{2}\La{)}\\[-1mm]
&&+\sum_{1\leq \rho\leq n}~~
-2i\langle\lnab{\mu}\be,\brho
\rangle\gf{\bbe}{\rho}{\bmu} 
-2i\langle\lnab{\mu}
\be,\rho\rangle\gf{\bbe}{\brho}{\bmu} 
-2i\langle\lnab{\mu}\bbe
,\brho\rangle\gf{\rho}{\be}{\bmu} \\[-1mm]
&&+\sum_{1\leq \rho\leq n}~~
2i\langle\lnab{\bbe}\mu,\brho\rangle\gf{\be}{\rho}{\bmu} 
-2i\langle\lnab{\mu}\bbe,\rho
\rangle\gf{\brho}{\be}{\bmu} 
+2i\langle\lnab{\bbe}
\mu,\rho\rangle\gf{\brho}{\be}{\bmu} \\[-1mm]
&&+\sum_{1\leq \rho\leq n}~~2i\langle\lnab{\bmu}\be
,\brho\rangle\gf{\bbe}{\rho}{\mu} +
2i\langle\lnab{\bmu}\be,\rho
\rangle\gf{\bbe}{\brho}{\mu} +2i\langle\lnab{\bmu}
\bbe,\brho\rangle\gf{\rho}{\be}{\mu} \\[-1mm]
&&+\sum_{1\leq \rho\leq n}~~-2i\langle\lnab{\bbe}
\bmu,\brho\rangle\gf{\rho}{\be}{\mu} 
+2i\langle\lnab{\bmu}\bbe
,\rho\rangle\gf{\brho}{\be}{\mu} 
-2i\langle\lnab{\bbe}
\bmu,\rho\rangle\gf{\brho}{\be}{\mu} \\[-1mm]
&&+\sum_{1\leq \rho\leq n}~~2i\langle\lnab{\bbe}\bmu
,\brho\rangle\gf{\be}{\mu}{\rho} +2i\langle\lnab{\bbe}
\bmu,\rho\rangle\gf{\be}{\mu}{\brho} -2i\langle\lnab{\bbe}
\mu,\rho\rangle\gf{\be}{\bmu}{\brho} \\[-1mm]
&&+\sum_{1\leq \rho\leq n}~~-2i\langle\lnab{\bbe}
\mu,\brho\rangle\gf{\be}{\bmu}{\rho} +2i\langle\lnab{\be}
\mu,\brho\rangle\gf{\bbe}{\rho}{\bmu} - 2i\langle\lnab{\be}
\mu,\brho\rangle\gf{\bbe}{\bmu}{\rho}\\[-1mm]
&&+\sum_{1\leq \rho\leq n}~~2i\langle\lnab{\be}
\bmu,\rho\rangle\gf{\bbe}{\mu}{\brho} -
2i\langle\lnab{\be}
\bmu,\rho\rangle\gf{\bbe}{\brho}{\mu}\\[-1mm]
&&-2\sum_{1\leq \rho\leq n}(\cos\theta_{\mu}-\cos\theta_{\rho})\la{(}|\langle\lnab{\be}
\mu,\brho\rangle |^{2} + |\langle\lnab{\bbe}
\mu,\brho\rangle |^{2}\la{)}.
\end{eqnarray*}
\end{Lm}
\em Proof.~~\em We denote by $\lnab{X}\lnab{Y}dF$ the 
covariant derivative of $\lnab{Y}dF$
in $T^{*}M\otimes F^{-1}TN$, and  by $\overline{R}(X,Y)\xi$, 
the curvature tensor of this connection, namely 
$(\overline{R}(X,Y)\xi)(Z)$
 $= R^{N}(dF(X),dF(Y))\xi (Z) -\xi (R^{M}(X,Y)Z)$.
 From Lemma 5.3, for $p$ on a neighbourhood of $p_{0}$,\\[-5mm]
\[
d\tg_{\mu\bmu}(\be)=i\gdf{\be}{\mu}{\bmu}-i\gdf{\be}{\bmu}{\mu}+ 2\sum_{\rho}\La{(}\langle\lnab{\be}\mu,\brho
\rangle\tg_{\rho\bmu}+\langle\lnab{\be}\bmu,\rho
\rangle\tg_{\mu\brho}\La{)}.\\[-4mm]
\]
Then at $p_{0}$,\\[-8mm]
\begin{eqnarray}
\lefteqn{d\la{(}d\tg_{\mu\bmu}(\be)\la{)}(\bbe)=~~~~~~~}\nonumber\\
&=& ig\La{(}\lnab{\bbe}\la{(}\lnab{\be}dF(\mu)\la{)},
JdF(\bmu)\La{)} + ig\La{(}\lnab{\be}dF(\mu),
\lnab{\bbe}\la{(}JdF(\bmu)\la{)}\La{)}\nonumber\\[-1mm]
&&- ig\La{(}\lnab{\bbe}\la{(}\lnab{\be}dF(\bmu)\la{)},
JdF(\mu)\La{)} - ig\La{(}\lnab{\be}dF(\bmu),
\lnab{\bbe}\la{(}JdF(\mu)\la{)}\La{)}\nonumber\\[-1mm]
&&+2\sum_{\rho}
\La{(}\lnab{\bbe}\la{(}\langle\lnab{\be}
\mu,\brho\rangle \la{)}\tg_{\rho\bmu}+
\lnab{\bbe}\la{(}\langle\lnab{\be}
\bmu,\rho\rangle \la{)}\tg_{\mu\brho}\La{)}\nonumber\\[-1mm]
&&+\sum_{\rho}2\langle\lnab{\be}
\mu,\brho\rangle d \tg_{\rho\bmu}(\bbe) 
+2\langle\lnab{\be}
\bmu,\rho\rangle d \tg_{\mu\brho}(\bbe)\\[1mm]
&=& ig\La{(}\lnab{\bbe}\la{(}\lnab{\be}dF(\mu)\la{)},
JdF(\bmu)\La{)}+ig\La{(}\lnab{\be}dF(\mu),
J\lnab{\bbe}dF(\bmu)\La{)}\nonumber\\[-1mm]
&&+ig\La{(}\lnab{\be}dF(\mu),JdF(\lnab{\bbe}\bmu)\La{)}
-ig\La{(}\lnab{\bbe}\la{(}\lnab{\be}dF(\bmu)\la{)},
JdF(\mu)\La{)}\nonumber\\[-1mm]
&&-ig\La{(}\lnab{\be}dF(\bmu),J\lnab{\bbe}dF(\mu)\La{)}
-ig\La{(}\lnab{\be}dF(\bmu),
JdF(\lnab{\bbe}\mu)\La{)}\nonumber\\[-1mm]
&&+\cos\theta_{\mu}\La{(}\lnab{\bbe}\la{(}\langle\lnab{\be}
\mu,\bmu\rangle \la{)}+\lnab{\bbe}\la{(}\langle \mu,\lnab{\be}\bmu\rangle \la{)}\La{)}+(5.8)\nonumber\\[1mm]
&=& ig\La{(}\lnab{\bbe}\la{(}\lnab{\be}dF(\mu)\la{)},
JdF(\bmu)\La{)}\\
&&+ig\La{(}\lnab{\be}dF(\mu),
J\lnab{\bbe}dF(\bmu)\La{)}
+\sum_{\rho}2i\langle\lnab{\bbe}\bmu,\rho\rangle
\gf{\be}{\mu}{\brho}+2i\langle\lnab{\bbe}\bmu,\brho\rangle
\gf{\be}{\mu}{\rho}\nonumber\\[-3mm]
&&-ig\La{(}\lnab{\bbe}\la{(}\lnab{\be}dF(\bmu)\la{)},
JdF(\mu)\La{)}\\
&&-ig\La{(}\lnab{\be}dF(\bmu),J\lnab{\bbe}dF(\mu)\La{)}
+\sum_{\rho}-2i\langle\lnab{\bbe}\mu,\rho\rangle
\gf{\be}{\bmu}{\brho}-2i\langle\lnab{\bbe}\mu,\brho\rangle
\gf{\be}{\bmu}{\rho}\nonumber\\[-4mm]
&&+\cos\theta_{\mu}\La{(}\lnab{\bbe}\la{(}\langle\lnab{\be}
\mu,\bmu\rangle \la{)}+\lnab{\bbe}\la{(}\langle \mu,\lnab{\be}\bmu\rangle \la{)}\La{)}\\[-1mm]
&&+(5.8). \non \\[-8mm]\non
\end{eqnarray} 
The term (5.11) vanish because $\langle\lnab{\be}
\mu,\bmu\rangle =- \langle \mu,\lnab{\be}
\bmu\rangle$ on a neighbourhood of $p_{0}$. 
Minimality of $F$ implies \\[-8mm]
\begin{eqnarray*}
\lefteqn{\sum_{\be}~\lnab{\bbe}\La{(}\lnab{\be}dF(\mu)\La{)}=}\\[-2mm]
&=&\sum_{\be}~\lnab{\bbe}\La{(}\lnab{\mu}dF(\be)\La{)}
=\sum_{\be}~\lnab{\bbe}\lnab{\mu}dF(\be)+
\lnab{\mu}dF(\lnab{\bbe}\be)\\[-1mm]
&=&\sum_{\be}~\lnab{\mu}\lnab{\bbe}dF(\be) -\lnab{[\mu,\bbe]}dF(\be)+(\overline{R}(\mu,\bbe)dF)
(\be)+\lnab{\mu}dF(\lnab{\bbe}\be)\\[-1mm]
&=&\sum_{\be}~~\lnab{\mu}\la{(}\lnab{\bbe}dF(\be)\la{)}
-\lnab{\bbe}dF(\lnab{\mu}\be)
-\lnab{[\mu,\bbe]}dF(\be)\\[-3mm]
&&~~~+R^{N}(dF(\mu),dF(\bbe))dF(\be)
-dF(R^{M}(\mu,\bbe)\be) +\lnab{\mu}dF(\lnab{\bbe}\be)\\[2mm]
&=& \sum_{\be}~\sum_{\rho}-2\langle \lnab{\mu}\be,\brho\rangle
\lnab{\bbe}dF(\rho)+\sum_{\rho}-2\langle \lnab{\mu}\be,\rho\rangle
\lnab{\bbe}dF(\brho)\\[-2mm]
&&~~~-\sum_{\rho}(2\langle \lnab{\mu}\bbe,\brho\rangle
-2\langle \lnab{\bbe}\mu,\brho\rangle)\lnab{\rho}dF(Z_
{\be})\\[-2mm]
&&~~~-\sum_{\rho}(2\langle \lnab{\mu}\bbe,\rho\rangle
-2\langle \lnab{\bbe}\mu,\rho\rangle)\lnab{\brho}dF(Z_
{\be})\\[-2mm]
&&~~~+R^{N}(dF(\mu),dF(\bbe))dF(\be)
-dF(R^{M}(\mu,\bbe)\be) \\
&&~~~ +\sum_{\rho}2\langle \lnab{\bbe}\be,\brho\rangle
\lnab{\mu}dF(\rho)+\sum_{\rho}2\langle \lnab{\bbe}\be,\rho\rangle
\lnab{\mu}dF(\brho).\\[-9mm]
\end{eqnarray*}
Hence
\begin{eqnarray*}
\lefteqn{(5.9)=\sum_{\be}~~iR^{N}(dF(\mu),dF(\bbe),dF(\be),JdF(\bmu))
-\cos\theta_{\mu}R^{M}(\mu,\bbe,\be,\bmu)} \\[-2mm]
&&~~~~+\sum_{\be\rho}~~-2i\langle \lnab{\mu}\be,\brho\rangle
\gf{\bbe}{\rho}{\bmu}-2i\langle \lnab{\mu}\be,\rho\rangle
\gf{\bbe}{\brho}{\bmu}\\[-2mm]
&&~~~~+\sum_{\be\rho}2i(-\langle \lnab{\mu}\bbe,\brho\rangle
+\langle \lnab{\bbe}\mu,\brho\rangle)\gf{\rho}{\be}{\bmu}
+2i(-\langle \lnab{\mu}\bbe,\rho\rangle
+\langle \lnab{\bbe}\mu,\rho\rangle)
\gf{\brho}{\be}{\bmu}\\[-2mm]
&& ~~~~+\sum_{\be\rho}~~2i\langle \lnab{\bbe}\be,\brho\rangle
\gf{\mu}{\rho}{\bmu}+2i\langle \lnab{\bbe}\be,\rho\rangle
\gf{\mu}{\brho}{\bmu}.\\[-9mm]
\end{eqnarray*}
Similarly
\begin{eqnarray*}
\lefteqn{-(5.10) =\sum_{\be}~~iR^{N}(dF(\bmu),dF(\bbe),dF(\be),JdF(\mu))
+\cos\theta_{\mu}R^{M}(\bmu,\bbe,\be,\mu) }\\[-2mm]
&&~~~~~~+\sum_{\be\rho}-2i\langle \lnab{\bmu}\be,\brho\rangle
\gf{\bbe}{\rho}{\mu}
-2i\langle \lnab{\bmu}\be,\rho\rangle
\gf{\bbe}{\brho}{\mu}\\[-2mm]
&&~~~~~~+\sum_{\be\rho}~~2i(-\langle \lnab{\bmu}\bbe,\brho\rangle
+\langle \lnab{\bbe}\bmu,\brho\rangle)\gf{\rho}{\be}{\mu}
+2i(-\langle \lnab{\bmu}\bbe,\rho\rangle
+\langle \lnab{\bbe}\bmu,\rho\rangle)
\gf{\brho}{\be}{\mu}\\[-2mm]
&&~~~~~~ +\sum_{\be\rho}2i\langle \lnab{\bbe}\be,\brho\rangle
\gf{\bmu}{\rho}{\mu}+2i\langle \lnab{\bbe}\be,\rho\rangle
\gf{\bmu}{\brho}{\mu}.
\end{eqnarray*}
Using Bianchi identity,\\[-6mm]
\begin{eqnarray*}
\lefteqn{
iR^{N}(dF(\mu),dF(\bbe),dF(\be),JdF(\bmu))
-iR^{N}(dF(\bmu),dF(\bbe),dF(\be),JdF(\mu))=}\\
&=& -iR^{N}(dF(\be),dF(\mu),dF(\bbe),JdF(\bmu))
-iR^{N}(dF(\bbe),dF(\be),dF(\mu),JdF(\bmu))\\[-2mm]
&&-iR^{N}(dF(\bmu),dF(\bbe),dF(\be),JdF(\mu))\\
&=&iR^{N}(dF(\be),dF(\bbe),dF(\mu),JdF(\bmu))
+2Im\La{(}R^{N}(dF(\be),dF(\mu),dF(\bbe),JdF(\bmu))\La{)},
\end{eqnarray*}
and by Gauss equation, and minimality of $F$,\\[-6mm]
\begin{eqnarray*}
\lefteqn{\sum_{\be}-R^{M}(\mu,\bbe,\be,\bmu)-
R^{M}(\bmu,\bbe,\be,\mu)=}\\[-2mm]
&=&\sum_{\be}R^{M}(\be,\mu,\bbe,\bmu)+
R^{M}(\bbe,\be,\mu,\bmu)-
R^{M}(\bmu,\bbe,\be,\mu)\\[-2mm]
&=&\sum_{\be} -R^{M}(\be,\bbe,\mu,\bmu)   +2R^{M}(\be,\mu,\bbe,\bmu)\\[-2mm]
&=&\sum_{\be}~~-R^{N}(dF(\be),dF(\bbe),dF(\mu),dF(\bmu))
\\[-4mm]
&&~~~~~~-g\La{(}\lnab{\be}dF(\mu),\lnab{\bbe}dF(\bmu)\La{)}
+g\La{(}\lnab{\be}dF(\bmu),\lnab{\bbe}dF(\mu)\La{)}\\
&&~~~~~~+2R^{N}(dF(\be),dF(\mu),dF(\bbe),dF(\bmu))\\
&&~~~~~~+2g\La{(}\lnab{\be}dF(\bbe),\lnab{\mu}dF(\bmu)\La{)}
-2g\La{(}\lnab{\be}dF(\bmu),\lnab{\mu}dF(\bbe)\La{)}\\[2mm]
&=&\sum_{\be}~~-R^{N}(dF(\be),dF(\bbe),dF(\mu),dF(\bmu))
+2R^{N}(dF(\be),dF(\mu),dF(\bbe),dF(\bmu))\\[-4mm]
&&~~~~~~-g\La{(}\lnab{\be}dF(\mu),\lnab{\bbe}dF(\bmu)\La{)}
-g\La{(}\lnab{\be}dF(\bmu),\lnab{\mu}dF(\bbe)\La{)}.
\end{eqnarray*}
Note that $R^{N}(dF(\be),dF(\mu),dF(\bbe),dF(\bmu)) =
Im\La{(}iR^{N}(dF(\be),dF(\mu),dF(\bbe),dF(\bmu))\La{)}$, since it is real.
Therefore,\\[-6mm]
\begin{eqnarray}
\lefteqn{\sum_{\be} d\la{(}
d\tg_{\mu\bmu}(\be)\la{)}(\bbe)=~~~~~~~}\nonumber\\[-2mm]
&=&\!\!\! \sum_{\be}~~iR^{N}(dF(\be),dF(\bbe),dF(\mu),JdF(\bmu)
+i\cos\theta_{\mu}dF(\bmu))\nonumber\\[-3mm]
&&~~+2Im\La{(}R^{N}(dF(\be),dF(\mu),dF(\bbe),JdF(\bmu)
+i\cos\theta_{\mu}dF(\bmu))\La{)}\nonumber\\
&&~-\cos\theta_{\mu}\,g\La{(}\lnab{\be}dF(\mu),\lnab{\bbe}
dF(\bmu)\La{)}
-\cos\theta_{\mu}\,g\La{(}\lnab{\be}dF(\bmu),\lnab{\mu}
dF(\bbe)\La{)}\\
&&~~+\sum_{\rho}-2i\langle \lnab{\mu}\be,\brho\rangle
\gf{\bbe}{\rho}{\bmu}-2i\langle \lnab{\mu}\be,\rho\rangle
\gf{\bbe}{\brho}{\bmu}\nonumber\\[-1mm]
&&~~+\sum_{\rho}2i(-\langle \lnab{\mu}\bbe,\brho\rangle
+\langle \lnab{\bbe}\mu,\brho\rangle)
\gf{\rho}{\be}{\bmu}+2i(-\langle \lnab{\mu}\bbe,\rho\rangle
+\langle \lnab{\bbe}\mu,\rho\rangle)
\gf{\brho}{\be}{\bmu}\nonumber\\[-1mm]
&&~~ +\sum_{\rho}2i\langle \lnab{\bbe}\be,\brho\rangle
\gf{\mu}{\rho}{\bmu}+2i\langle \lnab{\bbe}\be,\rho\rangle
\gf{\mu}{\brho}{\bmu}\\[-1mm]
&&~~+\sum_{\rho}2i\langle \lnab{\bmu}\be,\brho\rangle
\gf{\bbe}{\rho}{\mu} +2i\langle \lnab{\bmu}\be,\rho\rangle
\gf{\bbe}{\brho}{\mu}\nonumber\\[-1mm]
&&~~+\sum_{\rho}2i(\langle \lnab{\bmu}\bbe,\brho\rangle
-\langle \lnab{\bbe}\bmu,\brho\rangle)\gf{\rho}{\be}{\mu}
+2(\langle \lnab{\bmu}\bbe,\rho\rangle
-\langle \lnab{\bbe}\bmu,\rho\rangle)
\gf{\brho}{\be}{\mu}\nonumber\\[-2mm]
&&~~ +\sum_{\rho}-2i\langle \lnab{\bbe}\be,\brho\rangle
\gf{\bmu}{\rho}{\mu}-2i\langle \lnab{\bbe}\be,\rho\rangle
\gf{\bmu}{\brho}{\mu}\\[-2mm]
&&~~+ig\La{(}\lnab{\be}dF(\mu),
J\lnab{\bbe}dF(\bmu)\La{)}
-ig\La{(}\lnab{\be}dF(\bmu),J\lnab{\bbe}dF(\mu)\La{)}\\
&&~~+\sum_{\rho}2i\langle\lnab{\bbe}\bmu,\rho\rangle
\gf{\be}{\mu}{\brho}+2i\langle\lnab{\bbe}\bmu,\brho\rangle
\gf{\be}{\mu}{\rho}\nonumber\\[-1mm]
&&~~\sum_{\rho}-2i\langle\lnab{\bbe}\mu,\rho\rangle
\gf{\be}{\bmu}{\brho}-2i\langle\lnab{\bbe}\mu,\brho\rangle
\gf{\be}{\bmu}{\rho} ~+~ (5.8).\nonumber
\end{eqnarray}
Using the unitary basis $\{\frac{\sqrt{2}}{\sin\theta_{\rho}}\Phi(\rho),
\frac{\sqrt{2}}{\sin\theta_{\rho}}\Phi(\brho)\}$ of the normal bundle,
and (2.1)\\[-9mm]
\begin{eqnarray*}
\lefteqn{(5.12)+(5.15)=}\\
&=&-\sum_{\be,\rho}\frac{2\cos\theta_{\mu}}{\sin^{2}\theta_{\rho}}\La{(}
|\gf{\be}{\mu}{\rho}|^{2}+|\gf{\be}{\mu}{\brho}|^{2}\La{)}
-\sum_{\be,\rho}\frac{2\cos\theta_{\mu}}{\sin^{2}\theta_{\rho}}\La{(}
|\gf{\be}{\bmu}{\rho}|^{2}+|\gf{\be}{\bmu}{\brho}|^{2}\La{)}\\[-2mm]
&&-\sum_{\be,\rho}\frac{2\cos\theta_{\rho}}{\sin^{2}\theta_{\rho}}\La{(}
|\gf{\be}{\mu}{\brho}|^{2}-|\gf{\be}{\mu}{\rho}|^{2}\La{)}
+\sum_{\be,\rho}\frac{2\cos\theta_{\rho}}{\sin^{2}\theta_{\rho}}\La{(}
|\gf{\be}{\bmu}{\brho}|^{2}-|\gf{\be}{\bmu}{\rho}|^{2}\La{)}\\
&=&2\sum_{\be,\rho}\frac{(\cos\theta_{\rho}-\cos\theta_{\mu})}
{\sin^{2}\theta_{\rho}}\;|\gf{\be}{\mu}{\rho}|^{2}
-2\sum_{\be,\rho}\frac{(\cos\theta_{\rho}+
\cos\theta_{\mu})}{\sin^{2}\theta_{\rho}}\;|\gf{\be}{\mu}{\brho}|^{2}\\[-2mm]
&&-2\sum_{\be,\rho}\frac{(\cos\theta_{\rho}+\cos\theta_{\mu})}
{\sin^{2}\theta_{\rho}}\;|\gf{\be}{\bmu}{\rho}|^{2}
+2\sum_{\be,\rho}\frac{(\cos\theta_{\rho}-
\cos\theta_{\mu})}{\sin^{2}\theta_{\rho}}\;|\gf{\be}{\bmu}{\brho}|^{2}.
\end{eqnarray*}
Applying lemma 5.3  and since  $\langle\lnab{Z}\mu,\bmu\rangle
+\langle\lnab{Z}\bmu,\mu\rangle=0$, we have
\begin{eqnarray*}
d\tg_{\mu\bmu}(\lnab{\bbe}\be)&=&
\sum_{\rho}2\langle \lnab{\bbe}\be,\brho\rangle d\tg_{\mu\bmu}(\rho)+\sum_{\rho}2\langle \lnab{\bbe}\be,
\rho\rangle d\tg_{\mu\bmu}(\brho)\\[-2mm]
&=& 2i\sum_{\rho}\LA{(}\langle \lnab{\bbe}\be,
\brho\rangle \gf{\rho}{\mu}{\bmu}
-\langle \lnab{\bbe}\be,\brho\rangle \gf{\rho}{\bmu}{\mu}
+\langle \lnab{\bbe}\be,
\rho\rangle \gf{\brho}{\mu}{\bmu}
-\langle \lnab{\bbe}\be,
\rho\rangle \gf{\brho}{\bmu}{\mu}\LA{)}\\[-1mm]
&=& (5.13)+(5.14).\\[-12mm]
\end{eqnarray*}
Finally\\[-9mm]
\begin{eqnarray*}
(5.8)&=&\sum_{\rho}2\langle \lnab{\be}\mu,\brho\rangle(
i\gf{\bbe}{\rho}{\bmu}-i\gf{\bbe}{\bmu}{\rho})-2\langle \lnab{\be}\mu,\brho\rangle
(\cos\theta_{\rho}-\cos\theta_{\mu})\langle \lnab{\bbe}\rho,
\bmu\rangle\\[-2mm]
&&+\sum_{\rho}2\langle \lnab{\be}\bmu,\rho\rangle(
i\gf{\bbe}{\mu}{\brho}-i\gf{\bbe}{\brho}{\mu}) -2\langle \lnab{\be}\bmu,\rho\rangle
(\cos\theta_{\mu}-\cos\theta_{\rho})\langle \lnab{\bbe}\mu,\brho\rangle\\
&=&\sum_{\rho}2i\langle \lnab{\be}\mu,\brho\rangle
\gf{\bbe}{\rho}{\bmu}
-2i\langle \lnab{\be}\mu,\brho\rangle\gf{\bbe}{\bmu}{\rho}+
2i\langle \lnab{\be}\bmu,\rho\rangle
\gf{\bbe}{\mu}{\brho}
-2i\langle \lnab{\be}\bmu,\rho\rangle
\gf{\bbe}{\brho}{\mu}\\[-3mm]
&&-2\sum_{\rho}(\cos\theta_{\mu}-\cos\theta_{\rho})\La{(}|
\langle \lnab{\be}\mu,\brho\rangle|^{2}+|
\langle \lnab{\bbe}\mu,\brho\rangle|^{2}\La{)}.\\[-8mm]
\end{eqnarray*}
These expressions lead  to the expression of the lemma.~~~~~~\qed\\[5mm]
Finally, we have
\begin{Pp}
If F is minimal without complex directions, then for each $0\leq k\leq 2n$  at each $p_{0}\in\Omega^{0}_{2k}$,
\begin{eqnarray*}
\triangle \kappa
&=&4i\sum_{\be} Ricci^{N}(JdF(\be),dF(\bbe))\\[-3mm]
&&+\sum_{\be,\mu}\!\frac{32}{\sin^{2}\theta_{\mu}} 
Im \La{(}R^{N}(dF(\be),dF(\mu),dF(\bbe),
JdF(\bmu)+i\cos\theta_{\mu}dF(\bmu))\La{)}\\[-2mm]
&&-\sum_{\be,\mu,\rho}\frac{64(\cos\theta_{\mu}
\!+\!\cos\theta_{\rho})}{\sin^{2}\theta_{\mu}\sin^{2}\theta_{\rho}}
Re\la{(}\gf{\be}{\mu}{\brho}\gf{\bbe}{\rho}{\bmu}\la{)}\\
&& +\sum_{\be,\mu,\rho}
\frac{32(\cos\theta_{\rho}-\cos\theta_{\mu})}
{\sin^{2}\theta_{\mu}\sin^{2}\theta_{\rho}}\;(|\gf{\be}{\mu}{\rho}|^{2}
+|\gf{\bbe}{\mu}{\rho}|^{2})\\
&&+\sum_{\be,\mu,\rho}\frac{32(\cos\theta_{\mu}+\cos\theta_{\rho})}
{\sin^{2}\theta_{\mu}}\,\LA{(}
|\langle\lnab{\be}\mu,\rho\rangle|^{2} +
|\langle\lnab{\bbe}\mu,\rho\rangle|^{2}\LA{)}.
\end{eqnarray*}
\end{Pp}
\em Proof. \em  From (5.7) and Lemma 5.4 we get\\[-8mm]
\begin{eqnarray}
\lefteqn{2\triangle \kappa=}\nonumber\\[-1mm]
=&&\!\!\!\!\!\!+\sum_{\be,\mu,\rho}\frac{64(\cos\theta_{\mu}+
\cos\theta_{\rho})}{\sin^{2}\theta_{\mu}\sin^{2}\theta_{\rho}} d\tg_{\mu\brho}(\bbe)d\tg_{\rho\bmu}(\be)\nonumber\\[-1mm]
&&\!\!\!\!\!\!+\!\sum_{\be,\mu}\!\frac{32i}{\sin^{2}\theta_{\mu}}
R^{N}(dF(\be),dF(\bbe),dF(\mu),
JdF(\bmu)+i\cos\theta_{\mu}dF(\bmu)) \nonumber\\[-2mm]
&&\!\!\!\!\!\!+\!\sum_{\be,\mu}\!\frac{64}{\sin^{2}\theta_{\mu}} 
Im \La{(}R^{N}(dF(\be),dF(\mu),dF(\bbe),
JdF(\bmu)+i\cos\theta_{\mu}dF(\bmu))\La{)} \nonumber\\[-1mm]
&&\!\!\!\!\!\! +\sum_{\be,\mu,\rho}\frac{64(\cos\theta_{\rho}-\cos\theta_{\mu})}
{\sin^{2}\theta_{\mu}\sin^{2}\theta_{\rho}}\;(|\gf{\be}{\mu}{\rho}|^{2}
+|\gf{\be}{\bmu}{\brho}|^{2})\nonumber\\[-1mm]
&&\!\!\!\!\!\! -\sum_{\be,\mu,\rho}\frac{64(\cos\theta_{\rho}+\cos\theta_{\mu})}
{\sin^{2}\theta_{\mu}\sin^{2}\theta_{\rho}}\;(|\gf{\be}{\mu}{\brho}|^{2}
+|\gf{\be}{\bmu}{\rho}|^{2})\nonumber\\[-1mm]
&&\!\!\!\!\!\! +\sum_{\be,\mu,\rho}-\frac{64i}{\sin^{2}\theta_{\mu}}
\langle\lnab{\mu}\be,\brho\rangle\gf{\bbe}{\rho}{\bmu} 
-\frac{64i}{\sin^2\theta_{\mu}}\langle\lnab{\mu}
\be,\rho\rangle\gf{\bbe}{\brho}{\bmu} -\frac{64i}{\sin^2\theta_{\mu}}
\langle\lnab{\mu}\bbe,\brho\rangle\gf{\rho}{\be}{\bmu} \\[-1mm]
&&\!\!\!\!\!\!+\sum_{\be,\mu,\rho}\frac{64i}
{\sin^2\theta_{\mu}}\langle\lnab{\bbe}
\mu,\brho\rangle\gf{\be}{\rho}{\bmu} -\frac{64i}{\sin^2\theta_{\mu}}
\langle\lnab{\mu}\bbe,\rho\rangle\gf{\brho}{\be}{\bmu} +\frac{64i}{\sin^2\theta_{\mu}}\langle\lnab{\bbe}
\mu,\rho\rangle\gf{\brho}{\be}{\bmu} \\[-1mm]
&&\!\!\!\!\!\!+\sum_{\be,\mu,\rho}\frac{64i}{\sin^2\theta_{\mu}}
\langle\lnab{\bmu}\be,\brho\rangle\gf{\bbe}{\rho}{\mu} +\frac{64i}{\sin^2\theta_{\mu}}
\langle\lnab{\bmu}\be,\rho\rangle\gf{\bbe}{\brho}{\mu} +\frac{64i}{\sin^2\theta_{\mu}}\langle\lnab{\bmu}
\bbe,\brho\rangle\gf{\rho}{\be}{\mu} \\[-1mm]
&&\!\!\!\!\!\!+\sum_{\be,\mu,\rho}-\frac{64i}
{\sin^2\theta_{\mu}}\langle\lnab{\bbe}
\bmu,\brho\rangle\gf{\rho}{\be}{\mu} +\frac{64i}{\sin^2\theta_{\mu}}
\langle\lnab{\bmu}\bbe,\rho\rangle\gf{\brho}{\be}{\mu} -\frac{64i}{\sin^2\theta_{\mu}}\langle\lnab{\bbe}
\bmu,\rho\rangle\gf{\brho}{\be}{\mu} \\[-1mm]
&&\!\!\!\!\!\!+\sum_{\be,\mu,\rho}\frac{64i}{\sin^2\theta_{\mu}}
\langle\lnab{\bbe}\bmu,\brho\rangle\gf{\be}{\mu}{\rho} +\frac{64i}{\sin^2\theta_{\mu}}\langle\lnab{\bbe}
\bmu,\rho\rangle\gf{\be}{\mu}{\brho} -\frac{64i}{\sin^2\theta_{\mu}}\langle\lnab{\bbe}
\mu,\rho\rangle\gf{\be}{\bmu}{\brho} \\[-1mm]
&&\!\!\!\!\!\!+\sum_{\be,\mu,\rho}-\frac{64i}{\sin^2\theta_{\mu}}
\langle\lnab{\bbe}\mu,\brho\rangle\gf{\be}{\bmu}{\rho} +\frac{64i}{\sin^2\theta_{\mu}}\langle\lnab{\be}
\mu,\brho\rangle\gf{\bbe}{\rho}{\bmu} -\frac{64i}{\sin^2\theta_{\mu}}\langle\lnab{\be}
\mu,\brho\rangle\gf{\bbe}{\bmu}{\rho}\\[-1mm]
&&+\sum_{\be,\mu,\rho}\frac{64i}{\sin^2\theta_{\mu}}\langle\lnab{\be}
\bmu,\rho\rangle\gf{\bbe}{\mu}{\brho} -\frac{64i}{\sin^2\theta_{\mu}}\langle\lnab{\be}
\bmu,\rho\rangle\gf{\bbe}{\brho}{\mu}\\[-1mm]
&&-\sum_{\be,\mu,\rho}\frac{64(\cos\theta_{\mu}-\cos\theta_{\rho})}
{\sin^{2}\theta_{\mu}}
\la{(}|\langle\lnab{\be}\mu,\brho\rangle |^{2}
 + |\langle\lnab{\bbe}\mu,\brho\rangle |^{2}\la{)}.\nonumber
\end{eqnarray}
Interchanging $\rho$ with $\be$ in the first term of (5.16) (that we named
by (5.16)(1), and similarly to other equations),  
we see that $ (5.16)(1)+(5.17)(2)=0$. 
Interchanging $\rho$ with $\be$ in (5.18)(1), we get
$(5.18)(1)+(5.19)(2)=0$.
In (5.16)(2),  $\langle \lnab{\mu}\be,\rho\rangle$ is skew-symmetric
on $\rho$ and $\be$, and $\gf{\bbe}{\brho}{\bmu}$  is symmetric on
$\rho$ and $\be$. Hence  $(5.16)(2)=0$.
 Similarly $(5.16)(3)=(5.18)(2)=(5.18)(3)=0$.
If we interchange  $\rho$ with $\mu$ in (5.17)(1),\\[-4mm]
\[(5.17)(1)+(5.20)(2)=-\sum_{\be,\mu,\rho}\frac{64i(\sin^{2}\theta_{\mu}-
\sin^{2}\theta_{\rho})}{\sin^{2}\theta_{\mu}\sin^{2}\theta_{\rho}}
\langle\lnab{\bbe}\bmu,\rho\rangle\gf{\be}{\mu}{\brho}.\\[-1mm]\]
Interchanging
 $\rho$ with $\mu$ in (5.17)(3), we get\\[-4mm]
\[(5.17)(3)+(5.20)(3)=-\sum_{\be,\mu,\rho}\frac{64i(\sin^{2}\theta_{\mu}+
\sin^{2}\theta_{\rho})}{\sin^{2}\theta_{\mu}\sin^{2}\theta_{\rho}}
\langle\lnab{\bbe}\mu,\rho\rangle\gf{\be}{\bmu}{\brho}.\\[-2mm]\]
Interchanging $\rho$ with $\mu$ in (5.19)(1), we get\\[-4mm]
\[(5.19)(1)+(5.20)(1)=\sum_{\be,\mu,\rho}\frac{64i(\sin^{2}\theta_{\mu}+
\sin^{2}\theta_{\rho})}{\sin^{2}\theta_{\mu}\sin^{2}\theta_{\rho}}
\langle\lnab{\bbe}\bmu,\brho\rangle\gf{\be}{\mu}{\rho}.\]
Interchanging  $\rho$ with $\mu$ in (5.19)(3), we get\\[-3mm]
\[(5.19)(3)+(5.21)(1)=\sum_{\be,\mu,\rho}\frac{64i(\sin^{2}\theta_{\mu}-
\sin^{2}\theta_{\rho})}{\sin^{2}\theta_{\mu}\sin^{2}\theta_{\rho}}
\langle\lnab{\bbe}\mu,\brho\rangle\gf{\be}{\bmu}{\rho}.\\[-2mm]\]
Interchanging  $\rho$ with $\mu$ in (5.21)(2),\\[-3mm]
\[(5.21)(2)+(5.22)(1)=\sum_{\be,\mu,\rho}\frac{64i(-\sin^{2}\theta_{\mu}+
\sin^{2}\theta_{\rho})}{\sin^{2}\theta_{\mu}\sin^{2}\theta_{\rho}}
\langle\lnab{\be}\bmu,\rho\rangle\gf{\bbe}{\mu}{\brho}.\\[-2mm]\]
Interchanging $\rho$ with $\mu$ in (5.22)(2), we obtain\\[-3mm]
\[(5.22)(2)+(5.21)(3)=\sum_{\be,\mu,\rho}\frac{64i(\sin^{2}\theta_{\mu}-
\sin^{2}\theta_{\rho})}{\sin^{2}\theta_{\mu}\sin^{2}\theta_{\rho}}
\langle\lnab{\be}\mu,\brho\rangle\gf{\bbe}{\bmu}{\rho}.\]\\[-9mm]
Therefore,\\[-9mm]
\begin{eqnarray}
\lefteqn{2\triangle \kappa=}\nonumber\\[-2mm]
&=&\sum_{\be,\mu,\rho}\frac{64(\cos\theta_{\mu}+
\cos\theta_{\rho})}{\sin^{2}\theta_{\mu}\sin^{2}\theta_{\rho}} d\tg_{\mu\brho}(\bbe)d\tg_{\rho\bmu}(\be)\\[-1mm]
&&+\sum_{\be,\mu}\!\frac{32i}{\sin^{2}\theta_{\mu}}
R^{N}(dF(\be),dF(\bbe),dF(\mu),
JdF(\bmu)+i\cos\theta_{\mu}dF(\bmu))\\[-2mm]
&&+\!\sum_{\be,\mu}\!\frac{64}{\sin^{2}\theta_{\mu}} 
Im \La{(}R^{N}(dF(\be),dF(\mu),dF(\bbe),
JdF(\bmu)+i\cos\theta_{\mu}dF(\bmu))\La{)} \nonumber\\[-1mm]
&& +\sum_{\be,\mu,\rho}\frac{64(\cos\theta_{\rho}-\cos\theta_{\mu})}
{\sin^{2}\theta_{\mu}\sin^{2}\theta_{\rho}}\;|\gf{\be}{\mu}{\rho}|^{2}
\\[-1mm]
&& -\sum_{\be,\mu,\rho}\frac{64(\cos\theta_{\rho}+\cos\theta_{\mu})}
{\sin^{2}\theta_{\mu}\sin^{2}\theta_{\rho}}\;|\gf{\be}{\mu}{\brho}|^{2}
\\[-1mm]
&& -\sum_{\be,\mu,\rho}\frac{64(\cos\theta_{\rho}+\cos\theta_{\mu})}
{\sin^{2}\theta_{\mu}\sin^{2}\theta_{\rho}}\;|\gf{\be}{\bmu}{\rho}|^{2}
\\[-1mm]
&& +\sum_{\be,\mu,\rho}\frac{64(\cos\theta_{\rho}-\cos\theta_{\mu})}
{\sin^{2}\theta_{\mu}\sin^{2}\theta_{\rho}}\;|\gf{\be}{\bmu}{\brho}|^{2}
\\[-1mm]
&&-\sum_{\be,\mu,\rho}\frac{64i(\sin^{2}\theta_{\mu}-
\sin^{2}\theta_{\rho})}{\sin^{2}\theta_{\mu}\sin^{2}\theta_{\rho}}
\langle\lnab{\bbe}\bmu,\rho\rangle\gf{\be}{\mu}{\brho} \\[-1mm]
&&-\sum_{\be,\mu,\rho}\frac{64i(\sin^{2}\theta_{\mu}+
\sin^{2}\theta_{\rho})}{\sin^{2}\theta_{\mu}\sin^{2}\theta_{\rho}}
\langle\lnab{\bbe}\mu,\rho\rangle\gf{\be}{\bmu}{\brho}\\[-1mm]
&&+\sum_{\be,\mu,\rho}\frac{64i(\sin^{2}\theta_{\mu}+
\sin^{2}\theta_{\rho})}{\sin^{2}\theta_{\mu}\sin^{2}\theta_{\rho}}
\langle\lnab{\bbe}\bmu,\brho\rangle\gf{\be}{\mu}{\rho}\\[-1mm]
&&+\sum_{\be,\mu,\rho}\frac{64i(\sin^{2}\theta_{\mu}-
\sin^{2}\theta_{\rho})}{\sin^{2}\theta_{\mu}\sin^{2}\theta_{\rho}}
\langle\lnab{\bbe}\mu,\brho\rangle\gf{\be}{\bmu}{\rho}\\[-1mm]
&&+\sum_{\be,\mu,\rho}\frac{64i(-\sin^{2}\theta_{\mu}+
\sin^{2}\theta_{\rho})}{\sin^{2}\theta_{\mu}\sin^{2}\theta_{\rho}}
\langle\lnab{\be}\bmu,\rho\rangle\gf{\bbe}{\mu}{\brho}\\[-1mm]
&&+\sum_{\be,\mu,\rho}\frac{64i(\sin^{2}\theta_{\mu}-
\sin^{2}\theta_{\rho})}{\sin^{2}\theta_{\mu}\sin^{2}\theta_{\rho}}
\langle\lnab{\be}\mu,\brho\rangle\gf{\bbe}{\bmu}{\rho}\\[-1mm]
&&-\sum_{\be,\mu,\rho}\frac{64(\cos\theta_{\mu}-\cos\theta_{\rho})}
{\sin^{2}\theta_{\mu}}
\la{(}|\langle\lnab{\be}\mu,\brho\rangle |^{2}
 + |\langle\lnab{\bbe}\mu,\brho\rangle |^{2}\la{)}.
\end{eqnarray}
By Lemma 5.3,\\[-6mm]
\begin{eqnarray}
(5.23)&=&\sum_{\be,\mu,\rho}~~
\frac{64(\cos\theta_{\mu}+\cos\theta_{\rho})}
{\sin^{2}\theta_{\mu}\sin^{2}\theta_{\rho}}\,
\cdot\LA{(}i\gf{\bbe}{\mu}{\brho}
-i\gf{\bbe}{\brho}{\mu}-(\cos\theta_{\mu}\!-\!\cos\theta_{\rho})
\langle\lnab{\bbe}\mu,\brho\rangle\LA{)}\!\cdot\nonumber\\[-4mm]
&&~~~~~~~~\cdot\LA{(}i\gf{\be}{\rho}{\bmu}
-i\gf{\be}{\bmu}{\rho}-(\cos\theta_{\rho}\!-\!\cos\theta_{\mu})
\langle\lnab{\be}\rho,\bmu\rangle\LA{)}\nonumber\\[1mm]
&=&-\sum_{\be,\mu,\rho}\frac{64(\cos\theta_{\mu}+\cos\theta_{\rho})}
{\sin^{2}\theta_{\mu}\sin^{2}\theta_{\rho}}\gf{\bbe}{\mu}{\brho}
\gf{\be}{\rho}{\bmu}\nonumber\\[-1mm]
&&+\sum_{\be,\mu,\rho}\frac{64(\cos\theta_{\mu}+\cos\theta_{\rho})}
{\sin^{2}\theta_{\mu}\sin^{2}\theta_{\rho}}|\gf{\be}{\bmu}{\rho}|^{2}\\[-1mm]
&&+\sum_{\be,\mu,\rho}\frac{64i(\cos^{2}\theta_{\mu}-\cos^{2}\theta_{\rho})}
{\sin^{2}\theta_{\mu}\sin^{2}\theta_{\rho}}\gf{\bbe}{\mu}{\brho}
\langle\lnab{\be}\rho,\bmu\rangle\\[-1mm]
&&+\sum_{\be,\mu,\rho}\frac{64(\cos\theta_{\mu}+\cos\theta_{\rho})}
{\sin^{2}\theta_{\mu}\sin^{2}\theta_{\rho}}|\gf{\be}{\rho}{\bmu}|^{2}\\[-1mm]
&&-\sum_{\be,\mu,\rho}\frac{64(\cos\theta_{\mu}+\cos\theta_{\rho})}
{\sin^{2}\theta_{\mu}\sin^{2}\theta_{\rho}}\gf{\be}{\bmu}{\rho}
\gf{\bbe}{\brho}{\mu}\nonumber\\[-1mm]
&&-\sum_{\be,\mu,\rho}\frac{64i(\cos^{2}\theta_{\mu}-\cos^{2}\theta_{\rho})}
{\sin^{2}\theta_{\mu}\sin^{2}\theta_{\rho}}
\langle\lnab{\be}\rho,\bmu\rangle\gf{\bbe}{\brho}{\mu}\\[-1mm]
&&-\sum_{\be,\mu,\rho}\frac{64i(\cos^{2}\theta_{\mu}-\cos^{2}\theta_{\rho})}
{\sin^{2}\theta_{\mu}\sin^{2}\theta_{\rho}}
\langle\lnab{\bbe}\mu,\brho\rangle\gf{\be}{\rho}{\bmu}\\[-1mm]
&&+\sum_{\be,\mu,\rho}\frac{64i(\cos^{2}\theta_{\mu}-\cos^{2}\theta_{\rho})}
{\sin^{2}\theta_{\mu}\sin^{2}\theta_{\rho}}
\langle\lnab{\bbe}\mu,\brho\rangle\gf{\be}{\bmu}{\rho}\\[-1mm]
&&+\sum_{\be,\mu,\rho}\frac{64(\cos^{2}\theta_{\mu}-\cos^{2}\theta_{\rho})}
{\sin^{2}\theta_{\mu}\sin^{2}\theta_{\rho}}(\cos\theta_{\rho}
-\cos\theta_{\mu})\langle \lnab{\bbe}\mu,\brho\rangle
\langle\lnab{\be}\rho,\bmu\rangle.
\end{eqnarray}
Immediately we have, $~(5.27)+(5.36)=(5.32)+(5.41)=(5.33)+(5.37)=0$, ~
and interchanging $\mu$ with $\rho$ in (5.26), (5.34) and in (5.40), we get,
$~ (5.26)+(5.38)=(5.29)+(5.40)=(5.34)+(5.39)=0$.
Note that\\[-5mm]
\[
\sum_{\mu,\rho}\frac{(\cos\theta_{\mu}-\cos\theta_{\rho})}
{\sin^{2}\theta_{\mu}}|\langle\lnab{\be}\mu,\brho\rangle |^{2}
=\sum_{\mu,\rho}\frac{(\cos\theta_{\rho}-\cos\theta_{\mu})}
{\sin^{2}\theta_{\rho}}|\langle\lnab{\bbe}\mu,\brho\rangle |^{2}.\\[-2mm]\]
 Hence $ (5.35)+(5.42)=0$. Then,
\begin{eqnarray}
2\triangle \kappa
&=&\!\!\!\!\sum_{\be,\mu}\!\frac{32i}{\sin^{2}\theta_{\mu}}
R^{N}(dF(\be),dF(\bbe),dF(\mu),
JdF(\bmu)+i\cos\theta_{\mu}dF(\bmu))\\[-2mm]
&&\!\!\!\!+\!\sum_{\be,\mu}\!\frac{64}{\sin^{2}\theta_{\mu}} 
Im \La{(}R^{N}(dF(\be),dF(\mu),dF(\bbe),
JdF(\bmu)+i\cos\theta_{\mu}dF(\bmu))\La{)} \nonumber\\[-2mm]
&&+\sum_{\be,\mu,\rho}-\frac{64(\cos\theta_{\mu}
\!+\!\cos\theta_{\rho})}{\sin^{2}\theta_{\mu}\sin^{2}\theta_{\rho}}
(\gf{\bbe}{\mu}{\brho}\gf{\be}{\rho}{\bmu}
+\gf{\be}{\bmu}{\rho}\gf{\bbe}{\brho}{\mu})\\[-1mm]
&& +\sum_{\be,\mu,\rho}
\frac{64(\cos\theta_{\rho}-\cos\theta_{\mu})}
{\sin^{2}\theta_{\mu}\sin^{2}\theta_{\rho}}\;(|\gf{\be}{\mu}{\rho}|^{2}
+|\gf{\bbe}{\mu}{\rho}|^{2})\non\\[-1mm]
&&\!\!\!\!-\sum_{\be,\mu,\rho}\frac{64i(\sin^{2}\theta_{\mu}+
\sin^{2}\theta_{\rho})}{\sin^{2}\theta_{\mu}\sin^{2}\theta_{\rho}}
\langle\lnab{\bbe}\mu,\rho\rangle\gf{\be}{\bmu}{\brho}\\[-1mm]
&&\!\!\!\!+\sum_{\be,\mu,\rho}\frac{64i(\sin^{2}\theta_{\mu}+
\sin^{2}\theta_{\rho})}{\sin^{2}\theta_{\mu}\sin^{2}\theta_{\rho}}
\langle\lnab{\bbe}\bmu,\brho\rangle\gf{\be}{\mu}{\rho}.
\end{eqnarray}
Using Lemma 5.3, and interchanging $\rho$ by $\mu$ when necessary,
\begin{eqnarray*}
\lefteqn{(5.45)+(5.46)=}\\[-1mm]
&\!\!\!\!\!\!\!\!=&\!\!\!\!\!\!\!\!\!\!\!
\sum_{\be,\mu,\rho}\!\!-\frac{64i}{\sin^{2}\theta_{\rho}}
\langle\lnab{\bbe}\mu,\rho\rangle\gf{\be}{\bmu}{\brho}
-\frac{64i}{\sin^{2}\theta_{\mu}}
\langle\lnab{\bbe}\mu,\rho\rangle\gf{\be}{\bmu}{\brho}
+\frac{64i}{\sin^{2}\theta_{\mu}}
\langle\lnab{\bbe}\bmu,\brho\rangle\gf{\be}{\mu}{\rho}
+\frac{64i}{\sin^{2}\theta_{\rho}}
\langle\lnab{\bbe}\bmu,\brho\rangle\gf{\be}{\mu}{\rho}\\[-1mm]
&\!\!\!\!\!\!\!\!=&\!\!\!\!\!\!\!\!
\sum_{\be,\mu,\rho}\frac{-64i}{\sin^{2}\theta_{\mu}}
\langle\lnab{\bbe}\mu,\rho\rangle\La{(}
\gf{\be}{\bmu}{\brho}-\gf{\be}{\brho}{\bmu}\La{)}
+\sum_{\be,\mu,\rho}\frac{64i}{\sin^{2}\theta_{\mu}}
\langle\lnab{\bbe}\bmu,\brho\rangle
\La{(}\gf{\be}{\mu}{\rho}-\gf{\be}{\rho}{\mu}\La{)}\\[-1mm]
&\!\!\!\!\!\!\!\!=&\!\!\!\!\!\!\!\!
\sum_{\be,\mu,\rho}\frac{64}{\sin^{2}\theta_{\mu}}
\langle\lnab{\bbe}\mu,\rho\rangle
(\cos\theta_{\mu}+\cos\theta_{\rho})
\langle\lnab{\be}\bmu,\brho\rangle
+\frac{64}{\sin^{2}\theta_{\mu}}
\langle\lnab{\bbe}\bmu,\brho\rangle
(\cos\theta_{\mu}+\cos\theta_{\rho})
\langle\lnab{\be}\mu,\rho\rangle\\[-1mm]
&=&\sum_{\be,\mu,\rho}\frac{64(\cos\theta_{\mu}+\cos\theta_{\rho})}
{\sin^{2}\theta_{\mu}}\,\LA{(}
|\langle\lnab{\be}\mu,\rho\rangle|^{2} +
|\langle\lnab{\bbe}\mu,\rho\rangle|^{2}\LA{)}.
\end{eqnarray*}
Obviously\\[-9mm]
\[(5.44)=
\sum_{\be,\mu,\rho}\frac{-128(\cos\theta_{\mu}
\!+\!\cos\theta_{\rho})}{\sin^{2}\theta_{\mu}\sin^{2}\theta_{\rho}}
Re\la{(}\gf{\be}{\mu}{\brho}\gf{\bbe}{\rho}{\bmu}\la{)}.
\]
 From (1.4), (2.1),  and  the $J$-invariance of $Ricci$,
$~(5.43)= 8i\sum_{\be} Ricci^{N}(JdF(\be),dF(\bbe))$,
and the expression of the Proposition follows.
\\[7mm]
{\Large \bf{Acknowledgments}}\\[3mm]
 We would like to thank very much Professor James Eells for helpful discussions and encouragement.\\[6mm]
{\Large \bf{References}}\\[2mm]
[B]~~{\small M.\ Berger, \em Sur les groupes d'holonomie des
vari\'{e}t\'{e}s \`{a} connexion affine et des 
vari\'{e}t\'{e}s riemannienes, \em Bull.\
 Soc.\ Math.\ France {\bf 83} (1955), 279-330. }\\[0mm]
[Ch-W]~~{\small S.S.\ Chern \& J.G.\ Wolfson, \em Minimal surfaces by moving frames, \em Amer.\  J.\  Math. {\bf 105} (1983), 59-83.}\\[0mm] 
[E-L]~~{\small J.\ Eells \& L.\ Lemaire, \em Selected topics
in harmonic maps, \em  C.B.M.S. Regional Conf.\ Series 
{\bf 50}, A.M.S. (1983).}\\[0mm]
[Mi-Mo]~~{\small M.J.\ Micallef \& J.D.\ Moore, \em
Minimal two-spheres and the topology of manifolds with positive
curvature on totally isotropic two-planes, \em Annals of Math. {\bf 127}
(1988), 199-227}\\[0mm]
[O]~~{\small Y.\ Ohnita, \em Minimal surfaces with constant curvature
and K\"{a}hler angle in complex space forms, \em
Tsukuba J.\ Math. {\bf 13} No1 (1989), 191-207.}\\[0mm] 
[O-V]~~{\small Y.\ Ohnita \& G.\ Valli \em Pluriharmonic maps into
 compact Lie groups and factorization into unitons, \em
Proc.\ London Math.\ Soc. {\bf 61} (1990), 546-570.}\\[0mm] 
[S-V]~~{\small I.\ Salavessa \& G.\ Valli, \em
Broadly-Pluriminimal Submanifolds of K\"{a}hler-Einstein Manifolds, \em
preprint submitted for publication.}\\[0mm]
[W]~~{\small J.\ G.\ Wolfson, \em Minimal
Surfaces in K\"{a}hler Surfaces and Ricci Curvature, \em J. Diff. Geom,
{\bf 29} (1989), 281--294.}\\[4mm] 
\end{document}